\documentstyle[12pt]{article}
\begin{document}

\title{ \Large {Characteristic Classes on Grassmann Manifolds}
 \footnotetext{
\mbox{}\quad  \ \ {\it Subject classification}. \  14M15, 55R10, 55U30, 57T15.\\
\mbox{}\quad  \ \ E-mail: \ jwzhou@suda.edu.cn and ttwomei@126.com}}
\author{Jianwei Zhou and Jin Shi}
\date{\small Department of Mathematics, Suzhou University, Suzhou
215006, P. R. China \\ Suzhou Senior Technical Institute, Suzhou
215007, P.R. China }

\maketitle
\begin{abstract}  In this paper,   we use characteristic classes of the canonical
vector bundles and the Poincar\' {e} dualality to study the
structure of the real homology and cohomology groups of oriented
Grassmann manifold $G(k, n)$. Show that for $k=2$ or $n\leq 8$,
the  cohomology groups $H^*(G(k,n),{\bf R})$ are generated by the
first Pontrjagin class, the Euler classes of the canonical vector
bundles. In these cases, the Poincar\' {e} dualality:
$H^q(G(k,n),{\bf R}) \to H_{k(n-k)-q}(G(k,n),{\bf R})$ can be
given explicitly.

\end{abstract}

\baselineskip 15pt
\parskip 3pt

\vskip 1cm

\noindent{\bf \S 1. Introduction } \vskip 0.3cm

Let $G(k, n)$ be the Grassmann manifold formed by all oriented
 $k$-dimensional subspaces of Euclidean space ${\bf R}^n$. For any $\pi\in
G(k, n)$, there are orthonormal vectors $e_1, \cdots, e_k$ such
that $\pi $ can be represented by $e_1\wedge\cdots\wedge e_k$. Thus
$G(k, n)$ becomes a submanifold of the  space $\bigwedge^k({\bf R}^n)$, thus
 we can use  moving frame to study the Grassmann manifolds.

There are two canonical vector bundles $E=E(k,n)$ and $F=F(k,n)$ over $G(k, n)$
with fibres generated by vectors of the
subspaces and the vectors orthogonal to the subspaces  respectively. Then we have Pontrjagin classes
$p_i(E)$ and $p_j(F)$ with relationship
$$(1 + p_1(E) + \cdots )(1 + p_1(F) + \cdots )=1.$$
If $k$ or $n-k$ is even number, we have Euler classes
$e(E)$ or $e(F)$.

The oriented Grassmann manifolds are classifying spaces for oriented vector bundles. For any
oriented vector bundle
$\tau\colon\; \xi \to M$ with fibre type ${{\bf R}}^k$,
there is a map $g\colon\; M\to G(k, n)$ such that $\xi$
is isomorphic to the induced bundle $g^*E(k, n)$.
If the maps $g_1,g_2\colon\; M\to G(k, n)$ are homotopy, the
induced bundles $g_1^*E(k, n)$ and $g_2^*E(k, n)$ are isomorphic.
 Then the characteristic classes
of the vector bundle $\tau$  are the pullback of the  characteristic classes of the vector bundle
$E(k,n)$.

In this paper, we use characteristic classes  and the Poincar\' {e} dualality to
study the real homology and cohomology groups of  oriented Grassmann manifolds.
The characteristic classes of the canonical vector bundles can be represented by
curvatures and are harmonic forms on the
Grassmann manifolds, see [7], [8], [13]. For $k=2$ or $n\leq 8$,
we show that the cohomology groups $H^*(G(k,n), {\bf R})$ are generated by the first Pontrjagin
class $p_1(E)$ and the Euler classes $e(E(k,n)), \ e(F(k,n))$ if $k$ or $n-k$ are even.
In these cases, the Poincar\' {e} duality: $H^q(G(k,n), {\bf R}) \to H_{k(n-k)-q}(G(k,n), {\bf R})$
can be given explicitly.

In \S 2, we compute volumes of some homogenous spaces which are
needed in later discussion. In \S 3, we study the Poincar\' {e} duality on oriented
compact Riemannian manifolds.  The results are

{\bf Theorem 3.1} \ Let $\varphi_1,\cdots,\varphi_k$ be a basis of the
cohomology group $H^q(M )$ represented by harmonic forms.  Then
$(\psi_1,\cdots,\psi_k)=(*\varphi_1,\cdots,*\varphi_k)A^{-1}$ is a basis of $
H^{n-q}(M )$,
where $a_{ij}=(\varphi_i,\varphi_j)$.  The
Poincar\' {e} duality $D\colon\; H^q(M )  \to  H_{n-q}(M)$ is given by
$$D(\varphi_i)=[T_i],$$
where  $[T_1],\cdots, [T_k]\in
H_{n-q}(M, {\bf
R})$ are the dual  of $\psi_1,\cdots,\psi_k$.
Furthermore, if $[S_1],\cdots, [S_k]$ are the dual basis of $\varphi_1,\cdots,\varphi_k$,
then $$D(\psi_i)=(-1)^{q(n-q)}[S_i].$$

 The Poincar\' {e} polynomials of
Grassmann manifold $G(k, n)$ for $k=2$ and $n\leq 8$
are listed in \S 3, which give  the structure of real homology and
cohomology groups of Grassmann manifold.
In \S 4, we show that the tangent space of Grassmann manifold is isomorphic to tensor
product of the canonical vector bundles.
Then we use the splitting principle of the characteristic class to
study the relationship among these vector bundles,
show that the characteristic classes of the tangent bundle on Grassmann manifold can be
represented by that of canonical vector bundles.

In \S 5,  we study
$G(2, N)$, the main result are the following

{\bf   Theorem 5.5} \ (1) When $2k+2 < N$, $[{\bf C}P^k]$ and
 $e^k(E(2, N))$ are the generators of $H_{2k}(G(2,N), {{\bf Z}})$ and
 $H^{2k}(G(2,N), {{\bf Z}})$ respectively;

(2) When $2k+2 > N$, $[G(2,k+2)]$ and
$\frac 12e^k(E(2, N))$
are the generators of $ H_{2k}(G(2,N), {{\bf Z}})$ and
$H^{2k}(G(2,N),{{\bf Z}})$ respectively;

(3)   When $2k+2 < N$, $D(e^k(E(2, N)))=[G(2,N-k)]$;
when $2k+2 > N$,  $D(\frac 12 e^k(E(2, N)))=(-1)^{n-k}[{\bf C}P^{n-k-2}]$;

(4) $[{\bf C}P^n], [\overline{{\bf C}P}^n]$
and $\frac 12 \{(-1)^ne^n(E(2, 2n+2))\pm e(F(2, 2n+2))\}$
are  generators of $ H_{2n}(G(2,2n+2),{{\bf Z}})$ and $H^{2n}(G(2,2n+2), {{\bf Z}})$ respectively. Furthermore,
$$D(\frac 12 [(-1)^ne^n(E(2, 2n+2))+ e(F(2, 2n+2))]) =[{\bf C}P^n],$$
$$D(\frac 12 [(-1)^ne^n(E(2, 2n+2))- e(F(2, 2n+2))]) =[\overline{{\bf C}P}^n].$$

In \S 6, we study the Grassmann manifold $G(3, 6)$, the main result are

{\bf  Theorem 6.1} \ (1) $\frac 12 p_1(E(3,6))\in  H^{4}(G(3,6),{{\bf Z}})$ is a generator and
its Poincar\' {e} dual  $[SLAG]$ is a generator of $H_5(G(3,6),{{\bf Z}})$;

(2) $\frac 4{3\pi} *p_1(E(3,6))\in  H^{5}(G(3,6),{{\bf Z}})$ is a generator and
its Poincar\' {e} dual  $[G(2,4)]$ is a generator of $H_4(G(3,6),{{\bf Z}})$.

In \S 7, \S 8 we study the Grassmann manifold $G(3, 7)$ and $G(3, 8)$, the main results are Theorem 7.5 and
8.4. In \S 9, we study  $G(4, 8)$, the main results are

{\bf Theorem 9.4} \ (1) $e(E),  e(F),  \frac 12(p_1(E) +e(E) -e(F))\in H^4(G(4,8), {\bf Z})$
are the generators, their dual generators  are
$[{\bf C}P^2], [*{\bf C}P^2], [G(2,4)]\in H_4(G(4,8), {\bf Z})$;

(2) $\frac 12e^3(E), \frac 12 e^3(F),
\frac 12 p_1^3(E)$ and $[G(4,7)],  [G(3,7)], [CAY]$ are the generators of  $H^{12}
(G(4,8), {\bf Z})$ and $H_{12}
(G(4,8), {\bf Z})$  respectively;

(3) The Poincar\' {e} duals of $e(E),  e(F),  \frac 12(p_1(E) +e(E) -e(F))$ are
$$[G(4,7)], \ [G(3,7)], \ [CAY]+[G(4,7)]-[G(3,7)]$$ respectively.

{\bf Theorem 9.5} \ The characteristic classes
$$\frac 12(e^2(F)+  p_1(E)e(F)),\ \frac 12(e^2(F)-  p_1(E)e(F)),$$
$$ \frac 12(e^2(E) + p_1(E)e(E)), \  \frac 12(e^2(E) - p_1(E)e(E))$$ are the
 generators of $H^{8}
(G(4,8), {\bf Z})$. Their Poincar\' {e} duals are
 $$[ASSOC], \ [\widetilde {ASSOC}], \ [*ASSOC], \ [* {\widetilde {ASSOC}}]$$ respectively.

As application, in \S 5 and \S 9, we consider the Gauss maps of submanifolds in Euclidean spaces.
For example, if $g\colon\; M \to G(4,8)$ is the Gauss map of an
immersion $f\colon\; M \to {{\bf R}}^8$ of a compact oriented $4$
dimensional manifold,
we have
$$g_*[M] =\frac 12\chi(M)[G(4,5)] + \lambda [G(1,5)] + \frac 32 Sign(M)[G(2,4)],$$
where $\lambda = \frac 12\int_{M} \ e(F(4,8))$ and $Sign(M)$ be the Signature of $M$. $\lambda=0$
if $f$ is an imbedding.

In \S 10 we use Gysin sequence to compute the cohomology of the homogenous space $ASSOC=G_2/SO(4)$,
which had been studied by Borel and Hirzebruch [4].

The structure of the cohomology groups of infinite Grassmann manifold $G(k, {\bf R}^\infty)$ are simple,
they are generated by Pontrjagin classes and the Euler class (if $k$ is even) of the canonical vector bundle
freely,
see [12], p.179.

\vskip 1cm \noindent{\bf \S 2. The volumes of homogenous spaces} \vskip
0.3cm

Let $e_1,  e_2, \cdots,  e_n$ be orthonormal frame fields on
${\bf R}^n$ such that $G(k, n)$ be generated by $e_1\wedge\cdots \wedge e_k$
locally. The vectors $e_1,  e_2, \cdots,  e_n$ can be viewed as functions on
Grassmann manifold. Let
$\mbox {d} e_A = \sum\limits_{B=1}^n \ \omega_A^Be_B,  \ \ \omega_A^B
 =\langle \mbox {d}e_A,e_B\rangle$ be $1$ forms on $G(k,n)$.
From $\mbox {d}^2 e_A =0$, we have
$\mbox {d}\omega_A^B= \sum\limits_{C=1}^n \ \omega_A^C\wedge \omega_C^B.$
 By $$\mbox {d}(e_1\wedge \cdots \wedge
e_k)=\sum\limits_{i=1}^k\sum\limits_{\alpha=k+1}^n \omega_i^\alpha
E_{i\alpha},$$
$$E_{i\alpha}=e_1\cdots e_{i-1}e_\alpha e_{i+1}\cdots
e_k, \ i=1, \cdots, k,   \ \alpha=k+1, \cdots,n, $$
we know
 $E_{i\alpha}$ form a basis of $T_{e_1\cdots e_k}G(k, n)$,  $\omega_i^\alpha $ be
their dual basis.
$$\mbox {d}s^2=\langle \mbox {d}(e_1\wedge \cdots \wedge
e_k),\mbox {d}(e_1\wedge \cdots \wedge
e_k)\rangle = \sum\limits_{i=1}^k
\sum\limits_{\alpha=k+1}^n (\omega_i^\alpha)^2
$$ is the induced metric on $G(k, n)$.  Differential
$E_{i\alpha}=e_1\cdots e_{i-1}e_\alpha e_{i+1}\cdots
e_k$, by Gauss equation,  we get the Riemannian connection $\nabla$ on $G(k, n)$,
$$\nabla E_{i\alpha}=\sum\limits_{j=1}^k \ \omega_i^{j}E_{j\alpha}
         +\sum\limits_{\beta=k+1}^n \ \omega_\alpha^\beta E_{i\beta}.  $$

Grassmann manifold $G(k, n)$ is oriented,
the orientation is given by the volume form
$$\omega_1^{k+1}\wedge\omega_2^{k+1}\wedge\cdots\wedge
\omega_k^{k+1}\wedge\cdots\wedge \omega_1^{n}\wedge\omega_2^{n}\wedge\cdots\wedge
\omega_k^{n}.$$

For later use we compute the volume of some homogenous spaces, such as real Grassmann manifold $G(k, n)$,
complex Grassmann manifold $G_{\bf C}(k, n)$ and $SLAG= SU(n)/SO(n)$. We first compute
the volume of special orthogonal group $SO( n)$.

Let $gl(n,  {{\bf R}})$ be the set of all $n\times n$ real matrices with the
inner product
$$\langle X, Y\rangle =\mbox{tr} \, (XY^t)=\sum\limits_{A,B} \ X_{AB}Y_{AB}, \ \
X=(X_{AB}), Y=(Y_{CD})\in gl(n,  {{\bf R}}).$$
Then $gl(n,  {{\bf R}})$ is an Euclidean space and
$SO(n)$ is a Riemannian submanifold of $gl(n,  {{\bf R}})$. Represent the elements of $SO(n)$ by
 $G=(e_1,\cdots,e_n)^t$,
where $e_A$  the $A$-th row of $G$. The vectors $e_1,\cdots,e_n$ can be viewed as functions of $SO(n)$,
then $\omega_A^B=\langle\mbox {d} e_A,  e_B\rangle=\mbox {d}e_A\cdot e_B^t$ are
$1$ forms on $SO(n)$,
$\omega_A^B + \omega^B_A=0$.
 Let
$E_{BC}$ be the matrix with  $1$ in the $B$-th row, $C$-th column,
the others being zero. We have
$$\mbox {d}GG^{-1}=(\omega_A^B)= \sum\limits_{A,B}  \ \omega_A^B E_{AB}, \ \ \mbox {d}G=\sum
\limits_{A,B}  \
\omega_A^BE_{AB}G=\sum\limits_{A<B}  \ \omega_A^B(E_{AB}-E_{BA})G.$$
Then $\{(E_{AB}-E_{BA})G\}$ is a basis of
$T_GSO(n)$ and
$$\mbox{d}s^2=\langle \mbox {d}G, \mbox {d}G\rangle = 2\sum\limits_{A<B} \omega_A^B\otimes
\omega_A^B$$ is Riemannian metric on $SO(n)$.

{\bf  Proposition 2.1} \ The volume of $SO( n)$  is
$$V(SO(n))= 2^{\frac 12(n-1)}V(S^{n-1})V(SO(n-1))=2^{\frac 14n(n-1)}V(S^{n-1})
\cdots V(S^1).$$

{\bf Proof} \ Let $\bar e_n=(0,\cdots,0,1)$ be a fixed vector. The
map $\tau(G)=\bar e_nG=e_n$ defines a fibre bundle $\tau\colon\;
SO(n)\to S^{n-1}$ with fibres $SO(n-1)$.  By $ \mbox {d}e_n=\sum \ \omega_n^Ae_A$,
$$\mbox {d}V_{S^{n-1}}=\omega^n_1\cdots\omega^n_{n-1}$$
is the volume element of $S^{n-1}$.
The volume element of $SO(n)$ can be
represented by
$$\mbox {d}V_{SO(n)}=(\sqrt 2)^{\frac 12n(n-1)}\prod_{A<B} \ \omega_A^B= 2^{\frac 12 (n-1)}
 (\sqrt 2)^{\frac 12(n-1)(n-2)}\prod_{A<B<n} \ \omega_A^B\cdot \tau^*\mbox {d}V_{S^{n-1}},$$
restricting $(\sqrt 2)^{\frac 12(n-1)(n-2)}\prod_{A<B<n} \ \omega_A^B$ on the fibres of $\tau$
are the volume elements of the fibres. Integration $\mbox {d}V_{SO(n)}$
along the fibre of $\tau$ first, then on $S^{n-1}$, these shows
$$V(SO(n))= 2^{\frac 12(n-1)}V(S^{n-1})V(SO(n-1)). \ \ \ \Box $$

As we know $V(S^{m})=\frac {2\pi^{\frac {m+1}2}}{\Gamma(\frac {m+1}2)}$,
$$V(S^{2n-1})=\frac {2\pi^n}{(n-1)!}, \ \ V(S^{2n})=\frac {2^{2n+1}n!\pi^n}{(2n)!}.$$

To compute the volume of $G(k,n)$,
we use principle bundle $SO(n) \to G(k,n)$ with the Lie group $SO(k)\times SO(n-k)$
as fibres.
The proof is similar to that of  Proposition 2.1.

{\bf  Proposition 2.2} \ The volume of Grassmann manifold $G(k,n)$ is
$$V(G(k,n))=\frac {V(SO(n))}{2^{\frac 12k(n-k)}V(SO(k))V(SO(n-k))}=
\frac {V(S^{n-1})\cdots V(S^{n-k})}{V(S^{k-1})\cdots V(S^{1})}.$$

By simple computation, we have
$$V(G(2,n+2))=\frac {2(2\pi)^n}{n!}, \ \ V(G(3,6))=\frac 23 \pi^5,$$ $$V(G(3,7))=\frac {16}{45} \pi^6, \ \
V(G(3,8))=\frac 2{45} \pi^8, \ \ V(G(4,8))=\frac 8{135} \pi^8.$$

Now  we compute the volume of complex Grassmann manifold $G_{\bf C}(k, n)$.
Let $J$ be the natural complex structure on ${\bf C}^n= {\bf R}^{2n}$ and
$s_1,\cdots,s_k$ be Hermitian orthonormal basis of $\pi \in G_{\bf C}(k, n)$.
Let $e_{2i-1}, \ e_{2i}=Je_{2i-1}\in  {\bf R}^{2n}$
 be the realization vectors of $s_i, \ \sqrt{-1} \, s_i$ respectively.
Then $e_1e_2\cdots e_{2k-1}e_{2k}\in G(2k,2n)$, and $G_{\bf C}(k, n)$ becomes
a submanifold of $G(2k,2n)$.

Let $U(n)=\{G\in gl(n,  {\bf C}) \ | \ G\overline G^t=I\} $ be the unitary group,
the Hermitian inner product of $X=(X_{AB}), Y=(Y_{CD})\in gl(n,  {\bf C})$ be
$$\langle X, Y\rangle =\mbox{tr} \, ( X\overline Y^t)=\sum\limits_{A,B} \ X_{AB}\overline  Y_{AB}.$$
Let $G=(s_1,\cdots,s_n)^t\in U(n)$ represented by the
 rows of $G$,
$\omega_A^B=\langle \mbox {d}s_A,  s_B\rangle=\mbox {d}s_A\cdot \bar s_B^t$ be
$1$ forms on $U(n)$.
Let $\omega_A^B = \varphi_A^B +\sqrt {-1}\, \psi_A^B$. From
$\omega_A^B + \overline\omega_B^A=0$ we have $\varphi_A^B+\varphi^A_B=0, \ \psi_A^B-\psi^A_B=0$.
Then
\begin{eqnarray*}  && \mbox {d}G = \sum\limits_{A, B} \ \omega_A^B E_{AB}G \\
&  & = \sum\limits_{A<B}  \ \varphi_A^B(E_{AB}-E_{BA})G+\sqrt {-1}\,
\{\sum\limits_{A< B}  \ \psi_A^B(E_{AB}+E_{BA})G+
\sum\limits_{A}  \ \psi_A^AE_{AA}G\}, \end{eqnarray*}
and
$$\mbox{d}s^2=\langle \mbox {d}G, \mbox {d}G\rangle = 2\sum\limits_{A<B} (\varphi_A^B\otimes
\varphi_A^B+\psi_A^B\otimes
\psi_A^B)+
\sum\limits_{A} \psi_A^A\otimes
\psi_A^A$$ is Riemannian metric on $U(n)$. The volume element is
$$\mbox{d}V_{U(n)}= 2^{\frac 12 n(n-1)}\psi_1^1\cdots\psi_n^n \prod_{A<B} \ \varphi_A^B\psi_A^B.$$

{\bf  Proposition 2.3} \ (1) The volume of $U( n)$  is
$$V(U(n))=2^{n-1}V(S^{2n-1})V(U(n-1))=2^{\frac 12n(n-1)}V(S^{2n-1})V(S^{2n-3})\cdots V(S^{1});$$

(2) As Riemannian submanifold of $G(2k,2n)$, the volume of  $G_{\bf C}(k,n)$ is
$$V(G_{\bf C}(k,n))=\frac {V(U(n))}{V(U(k))V(U(n-k))}.$$

{\bf Proof} \ Let $\bar e_n=(0,\cdots,0,1)$ be a fixed vector. The
map  $\tau(G)= \bar e_nG=s_n$ defines a fibre bundle $\tau\colon\;
U(n)\to S^{2n-1}$ with fibre type $U(n-1)$. From $\mbox {d}s_n=\sum \ \omega_n^As_A$ and $
\omega_n^A= \varphi_n^A +\sqrt {-1}\, \psi_n^A, \ \varphi_n^n=0$,
we have the volume element of $S^{2n-1}$,
$$\mbox {d}V_{S^{2n-1}}=\varphi^n_1\psi_1^n\cdots\varphi^n_{n-1}\psi^n_{n-1}\psi^n_{n}.$$
Then  the volume element of $U(n)$ can be
represented by
$$\mbox {d}V_{U(n)}= 2^{n-1}
\tau^*\mbox {d}V_{S^{2n-1}}\cdot\mbox {d}V_{U(n-1)}.$$
These prove (1).

As noted above, the map $[s_1\cdots s_k]\mapsto e_1e_2\cdots
e_{2k-1}e_{2k}$ gives an imbedding of $G_{\bf C}(k, n)$ in
$G(2k,2n)$.  From $\mbox{d}s_i = \sum\limits_{j=1}^k \
\omega_i^js_j + \sum\limits_{\alpha=k+1}^n
 \ \omega_i^\alpha s_\alpha, \ \omega_i^j = \varphi_i^j +\sqrt {-1}\, \psi_i^j, \
 \ \omega_i^\alpha = \varphi_i^\alpha +\sqrt {-1}\, \psi_i^\alpha$, we have
$$\mbox{d}e_{2i-1} = \sum\limits_{j=1}^k \ (\varphi_i^je_{2j-1} +\psi_i^je_{2j})
 + \sum\limits_{\alpha=k+1}^n
 \ (\varphi_i^\alpha e_{2\alpha-1}+ \psi_i^\alpha e_{2\alpha}), $$
$$\mbox{d}e_{2i} = \sum\limits_{j=1}^k \ (\varphi_i^je_{2j} -\psi_i^je_{2j-1})
 + \sum\limits_{\alpha=k+1}^n
 \ (\varphi_i^\alpha e_{2\alpha}- \psi_i^\alpha e_{2\alpha-1}). $$
Then
$$\mbox{d} (e_1e_2\cdots e_{2k-1}e_{2k})= \sum\limits_{i,\alpha} \ \varphi_i^\alpha
(E_{2i-1\, 2\alpha-1}+E_{2i\, 2\alpha})+ \sum\limits_{i,\alpha} \ \psi_i^\alpha
(E_{2i-1\, 2\alpha}-E_{2i\, 2\alpha-1}),$$
$$\mbox{d} V_{G_{\bf C}(k, n)}= 2^{k(n-k)}\varphi_1^{k+1}\psi_1^{k+1}\cdots
\varphi_k^{n}\psi_k^{n}.$$

 The rest is similar to that of Proposition 2.1. \ \ \ $\Box$

The symmetric space $SLAG= SU(n)/SO(n)$ can be imbedded in $G(n,
2n)$ as follows. Let $\bar e_{2i-1}, \bar e_{2i} =J\bar e_{2i-1},
\ i=1,\cdots,n, $ be a fixed orthonormal basis of ${\bf
C}^n={{\bf R}}^{2n}$, $M = \{G(\bar e_1\bar e_3\cdots\bar
e_{2n-1}) \ | \ G\in SU(n)\subset SO(2n)\}$ is diffeomorphic to
$SLAG= SU(n)/SO(n)$.

{\bf  Proposition 2.4} \ (1) The volume of special unitary group $SU(n)$ is
$$V(SU(n))=2^{n-1}\sqrt {\frac n{n-1}}V(S^{2n-1})V(SU(n-1));$$

(2) The volume of $M= SLAG$ is
$$V(M)=\frac {V(SU(n))}{V(SO(n))}.$$

{\bf  Proof} \ The proof is similar to  Proposition 2.1. Let
$G=(s_1,\cdots,s_n)^t\in SU(n)$, $\omega_A^B=\mbox {d}s_A\cdot
\bar s_B^t$. From $\det G=1$ we have
 $\sum\limits_{A=1}^n \ \omega_A^A=0$, then $
\psi_n^n=-\sum\limits_{B\not =n} \ \psi_B^B$. The Riemannian
metric on $SU(n)$ is
\begin{eqnarray*} \mbox{d}s^2 & = & 2\sum\limits_{A<B} (\varphi_A^B\otimes
\varphi_A^B+\psi_A^B\otimes \psi_A^B)+
\sum\limits_{B\neq n } \psi_B^B\otimes
\psi_B^B + \psi_n^n\otimes
\psi_n^n \\
&=& 2\sum\limits_{A<B} (\varphi_A^B\otimes
\varphi_A^B+\psi_A^B\otimes \psi_A^B) \\
&& \quad +
(\psi_1^1,\cdots,\psi_{n-1}^{n-1})\left(\begin{array}{ccccc} 2 & 1 & \cdots & 1 \\ 1 & 2 & \cdots &
1 \\ \vdots & \vdots & \ddots & \vdots
\\ 1 & 1 & \cdots & 2 \end{array}\right)
\left(\begin{array}{ccccc} \psi_1^1 \\ \vdots \\ \psi_{n-1}^{n-1}
\end{array}\right). \end{eqnarray*} Then
$$\mbox{d}V_{SU(n)}= 2^{\frac 12 n(n-1)}\sqrt n
\psi_1^1\cdots\psi_{n-1}^{n-1} \prod_{A<B} \ \varphi_A^B\psi_A^B.$$
The volume of special unitary group $SU(n)$ is
$$V(SU(n))=2^{n-1}\sqrt {\frac n{n-1}}V(S^{2n-1})V(SU(n-1)).$$

Let $e_{2A-1}, \ e_{2A}=Je_{2A-1}$
 be the realization vectors of $s_A, \ \sqrt{-1} \, s_A$ respectively.
$M=SLAG$ is generated by $G(\bar e_1\bar e_3\cdots \bar
e_{2n-1})=e_1 e_3\cdots e_{2n-1}$,
$$\mbox{d} (e_1 e_3\cdots e_{2n-1})=\sum \ \psi_B^B(E_{2B-1 \, 2B} -E_{2n-1 \, 2n}) +
\sum\limits_{A< B}\ \psi_A^B(E_{2A-1 \, 2B} + E_{2B-1 \, 2A}),$$
$$\mbox{d}s^2 = 2\sum\limits_{A<B} \psi_A^B\otimes
\psi_A^B + 2\sum\limits_{B\neq n } \psi_B^B\otimes \psi_B^B +
\sum\limits_{B\neq C< n } \psi_B^B\otimes \psi_C^C.$$ Then
$$\mbox{d}V_{M}= 2^{\frac 14n(n-1)}\sqrt n
\psi_1^1\cdots\psi_{n-1}^{n-1} \prod_{A<B} \ \psi_A^B.$$ Let
$\tau\colon\; SU(n)\to M$ be the projection, the fibre of $\tau$
at $ e_1 e_3\cdots  e_{2n-1}$ is the group $SO(n)$. Restricting
$\mbox{d}s_i = \sum\limits_{j=1}^k \ \omega_i^js_j +
\sum\limits_{\alpha=k+1}^n \ \omega_i^\alpha s_\alpha$ on the
fibre of $\tau$, we have $\omega_i^\alpha=0$ and $\psi_i^j=0$,
then $\mbox{d}V_{SO(n)}= 2^{\frac 14 n(n-1)}\prod_{A<B} \
\varphi_A^B$ is the volume element of the fibres. These completes
the proof. \ \ \ \ $\Box$

\vskip 1cm \noindent{\bf \S 3. The Poincar\'{e} duality} \vskip
0.3cm

Let $M$ be a compact oriented Riemannian manifold  and
$H_q(M)=H_q(M, {{\bf R}})$ its $q$-th singular homology group,
 $H^q(M )=H^q(M, {{\bf R}})$ be the $q$-th de Rham cohomology
group. For any $[\xi]\in H^q(M  )$ and $[z]=[\sum \ \lambda_i
\sigma_i]\in H_q(M )$, we can define
$$[\xi]([z]) =\int_z \ \xi=\sum \ \lambda_i\int_{\sigma_i} \ \xi=\sum \ \lambda_i\int_{\triangle^q} \ \sigma_i^*\xi,$$
where every singular simplex $\sigma_i\colon\; \triangle^q  \to M$ is differentiable.
If $[\xi]\in H^q(M, {{\bf Z}})$ and $[z]\in H_q(M, {{\bf Z}})$,
we have $[\xi]([z])\in {{\bf Z}}$.
By universal coefficients Theorem, we have
$$H^q(M, {{\bf R}} )\cong \mbox {Hom} \, (H_q(M, {{\bf R}} ), {{\bf R}}),$$
and
$$H^q(M, {{\bf Z}})\cong \mbox {Hom} \, (H_q(M, {{\bf Z}}), {{\bf Z}})\oplus
\mbox{Ext} \, (H_{q-1}(M, {{\bf Z}}), {{\bf Z}}).$$
On the other hand, we have Poincar\' {e} duality
$$D\colon\; H^q(M , {{\bf R}})  \to  H_{n-q}(M, {{\bf R}}), \ \  n=\dim M.$$
For any $[\xi]\in H^q(M ), \ D[\xi]\in H_{n-q}(M )$,
 we have
$$[\eta](D[\xi])=\int_{D[\xi]} \ \eta =\int_{M} \ \xi\wedge\eta$$
for any $[\eta]\in H^{n-q}(M )$.

In the following, we use harmonic forms  to represent
the Poincar\' {e} duality. Let $\varphi_1,\cdots,\varphi_k$ be
the basis of $ H^q(M )$ and $[T_i]=D(\varphi_i)$ be their
Poincar\' {e} duals. By Hodge Theorem, we can assume that
$\varphi_1,\cdots,\varphi_k$ are all the harmonic forms on $M$.
Then $*\varphi_1,\cdots,*\varphi_k$ are  the harmonic forms on $M$
and form a basis of $H^{n-q}(M )$. Let
$$a_{ij}=(\varphi_i,\varphi_j)=\int_M \ \langle\varphi_i,\varphi_j\rangle \, \mbox {d} V_M
=\int_M \ \varphi_i\wedge* \varphi_j$$ be the inner product of
differential forms $\varphi_i,\varphi_j$. Let
$\psi_1,\cdots,\psi_k$ be the dual basis of
$[T_1],\cdots,[T_k]$, also represented by harmonic forms. Assuming
$\psi_j=\sum \ *\varphi_ib_{ij}$, by Poincar\' {e} duality,
$$\delta_{ij}=\int_{T_i} \ \psi_j =\int_M \ \varphi_i\wedge \psi_j=
\int_M \ \sum \ \varphi_i\wedge *\varphi_lb_{lj}=\sum \
a_{il}b_{lj}.$$ This shows $(b_{ij})=(a_{ij})^{-1}=A^{-1}$, we have
$$(\psi_1,\cdots,\psi_k)=(*\varphi_1,\cdots,*\varphi_k)A^{-1}.$$

{\bf Theorem 3.1} \ Let $\varphi_1,\cdots,\varphi_k$ be a basis of the
cohomology group $H^q(M )$ represented by harmonic forms.  Then
$(\psi_1,\cdots,\psi_k)=(*\varphi_1,\cdots,*\varphi_k)A^{-1}$ is a basis of $
H^{n-q}(M )$,
where $a_{ij}=(\varphi_i,\varphi_j)$.  The
Poincar\' {e} duality $D\colon\; H^q(M )  \to  H_{n-q}(M)$ is given by
$$D(\varphi_i)=[T_i],$$
where  $[T_1],\cdots, [T_k]\in
H_{n-q}(M, {\bf
R})$ are the dual  of $\psi_1,\cdots,\psi_k$.
Furthermore, if $[S_1],\cdots, [S_k]$ are the dual basis of $\varphi_1,\cdots,\varphi_k$,
then $$D(\psi_i)=(-1)^{q(n-q)}[S_i].$$

{\bf Proof} \ The equations $D(\psi_i)=(-1)^{q(n-q)}[S_i]$ follow from
$**\varphi_i=(-1)^{q(n-q)}\varphi_i$ and $(\varphi_i,\varphi_j)=(*\varphi_i,*\varphi_j).$
 \ \ \ \ $\Box$

Theorem 3.1 can be applied to the
Poincar\' {e} duality $D\colon\; H^q(M, {\bf Z})  \to  H_{n-q}(M, {\bf Z})$
if  we ignore the torsion elements of $H^q(M, {\bf Z})$.

The $q$-th Betti number is the common dimension of the real
homology and coholomogy groups $H_q(G(k,n))$ and  $H^q(G(k,n) )$
(and is also the rank of $H_q(G(k,n), {\bf Z})$ and
$H^q(G(k,n), {\bf Z})$). The Poincar\' {e}  polynomials, with
the Betti numbers as coefficients are given by the following Table
(see[7], [8]).
\vskip 0.2cm
$$\begin{tabular}{|c|c|c|c|}\hline \ \
Grassmannian \ \  & Poincar\' {e}  polynomial  \\ \hline
\ \ $G(1, n+1)$ \ \ &    $1+t^n$    \\   \hline
        $ G(2, 2n+1)$ &  $1+t^2+t^4+\cdots +t^{4n-2}$              \\  \hline
 \ \    $ G(2, 2n+2)$  \ \ & \ \  $(1+t^{2n})(1+t^2+\cdots+t^{2n})$ \ \              \\  \hline
$ G(3, 6)$ &  $(1+t^4)(1+t^{5})$              \\  \hline
$ G(3, 7)$ &  $(1+t^4+t^8)(1+t^{4})$              \\  \hline
 $ G(3, 8)$ &  $(1+t^4+t^8)(1+t^{7})$              \\  \hline
$ G(4, 8)$ &  $(1+t^4+t^8)(1+t^{4})^2$              \\  \hline \end{tabular}$$

\vskip 1cm \noindent{\bf \S 4. The vector bundles on $G(k,n)$ } \vskip
0.3cm

Let $\tau_1\colon\;  E(k,n) \to G(k,n)$ be the canonical vector bundle on Grassmann manifold
$G(k, n)$, the fibre over
$\pi\in G(k,n)$ be the vectors of $\pi$. $E=E(k,n) $ is a Riemannian vector bundle
 with the induced metric.
Let $e_1,\cdots,e_k,e_{k+1},\cdots,e_n$ be orthonormal frame fields on ${\bf R}^n$,
$G(k, n)$ is locally generated by $e_1\cdots e_k=e_1\wedge \cdots\wedge e_k$.
Then $e_1,\cdots,e_k$ are local orthonormal  sections of the vector
bundle $\tau_1$. From
$\mbox {d} e_i
=\sum\limits _{j=1}^k \ \omega_i^je_j + \sum\limits _{\alpha=k+1}^n \ \omega_i^\alpha e_\alpha$,
we know that $\nabla e_i
=\sum\limits _{j=1}^k \ \omega_i^je_j$ defines a  Riemannian connection on $\tau_1$.
From $\nabla^2 e_i
=\sum\limits _{j=1}^k \ (\mbox{d}\omega_i^j-\sum\limits _{l=1}^k \ \omega_i^l\wedge\omega_l^j)e_j$, we have
curvature forms $$\Omega_i^j=\mbox{d}\omega_i^j-\sum\limits_{l=1}^k \
\omega_i^l\wedge\omega_l^j=\sum_{\alpha=k+1}^n \ \omega_i^\alpha\wedge\omega_\alpha^j.$$
The total Pontrjagin classes of the vector bundle $\tau_1\colon\;  E \to G(k,n)$ are
$$p(E)=1+p_1(E)+p_2(E)+\cdots =\det (I+\frac 1{2\pi} (\Omega_i^j)).$$
If $k$ is even, we have Euler class,
$$e(E)
= \frac {(-1)^{\frac k2}}{(4\pi)^{\frac k2}(\frac k2)!}
\sum\limits_{i_1,\cdots, i_{k}}  \  \varepsilon(i_1i_2\cdots i_{k})
\Omega_{i_1i_2}\Omega_{i_3i_4}\cdots \Omega_{i_{k-1}i_{k}}. $$

Similarly, we can define vector bundle $\tau_2\colon\;  F=F(k,n) \to G(k,n)$ on Grassmann manifold
$G(k, n)$, the fibre over
$e_1\cdots e_k\in G(k,n)$ is  the vectors  orthogonal to $e_1,\cdots, e_k$.
Then $e_{k+1},\cdots,e_{n}$ are local orthonormal sections of $F$. From
$\mbox {d} e_\alpha
=\sum\limits _{\beta=k+1}^n \ \omega_\alpha^\beta e_\beta
+ \sum\limits _{i=1}^k \ \omega_\alpha^ie_i$, we have  Riemannian connection
$\nabla e_\alpha
=\sum\limits _{\beta={k+1}}^n \ \omega_\alpha^\beta e_\beta.$
The curvature forms  are given by
$$\nabla^2 e_\alpha
=\sum\limits _{j=1}^k \ \Omega_\alpha^\beta e_\beta, \ \
\Omega_\alpha^\beta =\sum_{i=1}^k \ \omega_\alpha^i\wedge\omega_i^\beta.$$
The total Pontrjagin classes of the vector bundle $\tau_2\colon\;  F \to G(k,n)$ are
$$p(F)=1+p_1(F)+p_2(F)+\cdots =\det (I+\frac 1{2\pi} (\Omega_\alpha^\beta)).$$

The direct sum $E(k,n)\oplus F(k,n)= G(k,n)\times {\bf R}^n$ is  trivial,
we have
$$(1+p_1(E)+p_2(E)+\cdots )\cdot (1+p_1(F)+p_2(F)+\cdots )=1.$$
Then Pontrjagin classes of $F$ are determined by that of $E$. For example,
we have
$$p_1(F)=-p_1(E), \ \ p_2(F)= p_1^2(E)-p_2(E), \ \ p_3(F)= - p_1^3(E)+2p_1(E)p_2(E)-p_3(E).$$
Let $*\colon\; \bigwedge^k( {\bf R}^n) \to \bigwedge^{n-k}({\bf R}^n)$ be the star operator,
$*G(k,n)=G(n-k,n)$ and the canonical vector bundles  $E(k,n), F(k,n)$ are interchanged under the map $*$.

{\bf Proposition 4.1} \  The tangent space $TG(k,n)$ of Grassmann manifold is isomorphic to tensor
product  $E(k,n)\otimes F(k,n)$.
If $k(n-k)$ is even, we have
$$e(G(k,n))=e(E(k,n)\otimes F(k,n)).$$

{\bf Proof} \ Let $e_1,  e_2,  \cdots,  e_n$ be an oriented orthonormal basis of
${\bf R}^n$,  the fibre of $E(k,n)$ over $x=e_1\wedge \cdots\wedge e_k\in G(k,n)$ is generated by
  $e_1,\cdots, e_k$ and the fibre of $F(k,n)$ over $x$ is generated by
  $e_{k+1},\cdots, e_n$.
On the other hand, the tangent space $T_{x}G(k,n)$ is generated by
$E_{i \, \alpha}=
e_1\wedge\cdots\wedge e_{i-1} \wedge e_\alpha\wedge e_{i+1}\cdots\wedge e_k$.
It is easy to see that the map $E_{i \, \alpha}\mapsto e_i\otimes e_\alpha$ gives an isomorphism
from tangent bundle $TG(k,n)$ to tensor product $E\otimes F$.

It is easy to see that the isomorphism $TG(k, n) \to E(k,n)\otimes F(k,n)$  preserves the connections
on $TG(k,n)$ and $E\otimes F$ respectively, where the connection on $E\otimes F$ is
$$\nabla \, (e_i\otimes e_\alpha) = \sum\limits_{j=1}^k  \  \omega^j_ie_j\otimes e_\alpha +
 \sum\limits_{\alpha=k+1}^n
  \  \omega^\beta_\alpha e_i\otimes e_\beta.  \ \ \ \ \Box$$

In the following we use the splitting principle of the characteristic class
 to study the relationship among these vector bundles.
First we study the oriented Grassmann manifold $G(2k,2n)$.
Assuming $s_{1},\cdots,s_{2k}$ be the orthonormal sections of  vector bundle $E(2k, 2n)$ such that
the curvature of Riemannian connection has the form
 $$\frac 1{2\pi}\nabla^2
\left (
  \begin{array}{ccccc}
      s_{1}\\
      s_{2}\\
      \vdots\\
      s_{2k-1}\\
      s_{2k}
  \end{array}
  \right )
  =
   \left (
  \begin{array}{ccccc}
    0 & -x_1\\
   x_1  & 0\\
   & & \ddots\\
  & & &  0 & -x_k\\
  & & &  x_k  & 0\\
  \end{array}
  \right )
  \left (
  \begin{array}{ccccc}
      s_{1}\\
      s_{2}\\
       \vdots\\
     s_{2k-1}\\
      s_{2k}
  \end{array}
  \right ). $$
The total Pontrjagin classes and the Euler class of $E=E(2k,2n)$ are
$$p(E)=  \prod\limits_{i=1}^k \ (1+x_i^2), \ \ \
e(E)=  x_1\cdots x_k.$$

Similarly assuming $t_{2k+1},t_{2k+2},\cdots,t_{2n}$
be the orthonormal sections of  vector bundle $F(2k, 2n)$, \
the curvature of Riemannian connection   has the form
$$\frac 1{2\pi}\nabla^2
\left (
  \begin{array}{ccccc}
      t_{2k+1}\\
      t_{2k+2}\\
      \vdots\\
      t_{2n-1}\\
      t_{2n}
  \end{array}
  \right )
  =
   \left (
  \begin{array}{ccccc}
    0 & -y_{k+1}\\
   y_{k+1}  & 0\\
   & & \ddots\\
  & & &  0 & -y_n\\
  & & &  y_n  & 0\\
  \end{array}
  \right )
  \left (
  \begin{array}{ccccc}
      t_{2k+1}\\
      t_{2k+2}\\
      \vdots\\
      t_{2n-1}\\
      t_{2n}
  \end{array}
  \right ). $$
The total Pontrjagin classes and the Euler class of $F=F(2k,2n)$ are
$$p(F)=  \prod\limits_{\alpha=k+1}^n \ (1+ y_\alpha^2),  \ \ \ e(F)=  y_{k+1}\cdots y_n.$$

 $ s_{2i-1}\otimes  t_{2\alpha-1},
s_{2i}\otimes  t_{2\alpha-1},
s_{2i-1}\otimes  t_{2\alpha},
 s_{2i}\otimes  t_{2\alpha}$ are the local orthonormal sections of
vector bundle $E\otimes F\cong TG(2k,2n)$.
The curvature of Riemannian connection on  $E\otimes F$ is given by
 $$\frac 1{2\pi}\nabla^2
 \left (
 \begin{array}{ccccccccc}
   s_{2i-1}\otimes  t_{2\alpha-1}\\
   s_{2i}\otimes  t_{2\alpha-1}\\
   s_{2i-1}\otimes  t_{2\alpha}\\
   s_{2i}\otimes  t_{2\alpha}
 \end{array}
 \right )
 =  \left (
 \begin{array}{ccccccccc}
    0 & -x_i & -y_\alpha & 0\\
    x_i & 0 & 0 &  -y_\alpha \\
   y_\alpha & 0 & 0 & -x_i\\
   0 & y_\alpha & x_i & 0
 \end{array}
 \right )
 \left (
 \begin{array}{ccccccccc}
 s_{2i-1}\otimes  t_{2\alpha-1}\\
   s_{2i}\otimes  t_{2\alpha-1}\\
   s_{2i-1}\otimes  t_{2\alpha}\\
   s_{2i}\otimes  t_{2\alpha}
 \end{array}
 \right ). $$
Then we have

{\bf Lemma 4.2} \  (1) $e(TG(2k, 2n)) =e(E\otimes F)=\prod_{i,\alpha} \
  (x_i^2-y_\alpha^2);$

(2) $p(TG(2k, 2n)) = p(E\otimes F)=\prod_{i,\alpha} \
  (1+ 2(x_i^2+y_\alpha^2)+ (x_i^2-y_\alpha^2)^2).$

The other cases can be discussed similarly.

By simple computation, we have
$$p_1(TG(2k, 2n)) = (2n-2k)p_1(E) +2kp_1(F)= 2(n-2k)p_1(E).$$
In particular,  $p_1(TG(2k, 4k))=0$.

In next section, we shall show
$$e(TG(2, 2n+2)) =(n+1)e^{2n}(E(2, 2n+2)),$$ $$e(TG(2, 2n+3)) =(n+1)e^{2n+1}(E(2, 2n+3)).$$
We can also show
$$e(TG(3, 7)) =3e^{3}(F(3,7)), \ \ e(TG(4, 8)) =6e^{4}(E(4,8))=6e^{4}(F(4,8)).$$

\vskip 1cm \noindent{\bf \S 5. The cases of $G(2,N)$ } \vskip
0.3cm

In this section, we study the real homology and cohomology of Grassmann manifold $G(2,N)$.

As is well-known, the oriented Grassmann manifold $G(2, N)$ is Kaehler manifold and can be imbedded in complex
projective space. Here we give a new proof.
Let $e_1,e_2$ be oriented orthonormal basis of $\pi\in G(2,N)$, \ $e_1\mapsto e_2,\ e_2\mapsto -e_1$
 defines an almost complex structure
$$J\colon\;  T_{\pi}G(2,N) \to T_{\pi}G(2,N), $$
$$E_{1 \, \alpha}=e_\alpha\wedge e_2 \mapsto -e_\alpha\wedge e_1= E_{2 \, \alpha}, \
E_{2 \, \alpha}=e_1\wedge e_\alpha \mapsto e_2\wedge e_\alpha= -E_{1 \,\alpha}.$$
It is easy to see that $ J$ is well defined and preserves the metric on $G(2,N)$.

 {\bf Proposition 5.1} \  $G(2, N)$ is a Kaehler manifold with complex structure $J$.

{\bf Proof } \ Let $\nabla$ be the Riemannian connection on $TG(2,N)$ defined above. We have
$$(\nabla  J)E_{i \, \alpha}=\nabla ( JE_{i  \, \alpha})-J(\nabla E_{i \, \alpha})
=0, \ \ i=1,2.$$
Hence, $\nabla  J=0$,
$J$ is a complex structure and $G(2, N)$ is a  Kaehler manifold.  \ \ \ \ $\Box$

The Euler classes of canonical vector bundles $E=E(2,2n+2)$ and $F=F(2,2n+2)$ can be represented by
$$e(E) =\frac 1{2\pi}\sum\limits_{\alpha=3}^{2n+2} \ \omega_1^\alpha\wedge\omega_2^\alpha,$$
$$e(F) =\frac {(-1)^n}{(4\pi)^{n}n!}\sum\limits \ \varepsilon(\alpha_1\alpha_2\cdots\alpha_{2n})
\Omega_{\alpha_1\alpha_2}\wedge\cdots\wedge \Omega_{\alpha_{2n-1}\alpha_{2n}}.$$

$G(2, k+2)$ is a submanifold of $G(2,2n+2)$ whose elements
 contained in a fixed $k+2$ dimensional subspace of ${{\bf R}}^{2n+2}$,
$i\colon\; G(2, k+2)\to G(2,2n+2)$ the inclusion. Then,
$E(2,k+2)=i^*E(2, 2n+2)$ and $e(E(2,k+2))=i^*e(E(2, 2n+2))$.
Let $G(1, 2n+1)$ be a submanifold of $G(2,2n+2)$ with elements  $e_1
\wedge \bar e_2, \ \bar e_2=(0,\cdots,0,1)$, $j\colon\; G(1, 2n+1)\to G(2,2n+2)$ be the inclusion.

{\bf  Theorem 5.2} \ For Grassmann manifold $G(2,2n+2)$, we have

(1) $ p_q(F)= (-1)^qp_1^q(E)= (-1)^qe^{2q}(E), \ \
q=1,\cdots,n;$

(2) The Pontrjagin classes and Euler class of tangent bundle $TG(2,2n+2)$  can be represented by the
Euler class of $E$,
$$p_1(G(2, 2n+2))=2(n-1)e^2(E), \ \  p_2(G(2, 2n+2))=(2n^2-5n +9)e^4(E), \ \cdots,$$
 $$e(G(2, 2n+2)) =(n+1)e^{2n}(E);$$

(3) $\int_{G(1,2n+1)} \ e(F)=2, \ \ \int_{G(2,k+2)} \ e^k(E)=2, \ \ \ k=1,\cdots,2n.$

{\bf Proof} \ For Grassmann manifold $G(2,2n+2)$, we have
$p_1(E)=e^2(E), \ p_{n}(F)=e^2(F)$. From $(1+p_1(E))\cdot
(1+p_1(F)+p_2(F)+\cdots+p_n(F) )=1$, we have
$$1+p_1(F)+p_2(F)+\cdots+p_n(F)= \frac 1{1+p_1(E)}= 1+ \sum\limits_{q=1}^n \ (-1)^{q}p_1^q(E).$$
Hence,
$p_q(F)= (-1)^{q}p_1^q(E)= (-1)^qe^{2q}(E)$, \
$p_{n}(F)=e^2(F)=(-1)^ne^{2n}(E)$. These proves (1).

By Lemma 4.2, note that $x_1=e(E)$, the Euler class of $G(2,2n+2)$ is
\begin{eqnarray*} e(TG(2, 2n+2)) & = & (x_1^2-y_2^2)(x_1^2-y_3^2)\cdots (x_1^2-y_{n+1}^2) \\
& = & x_1^{2n} -x_1^{2n-2}p_1(F) +x_1^{2n-4}p_2(F)-\cdots +(-1)^np_{n}(F) \\
& = & (n+1)e^{2n}(E). \end{eqnarray*}
By Gauss-Bonnet formula, we have
$$\chi(G(2,2n+2))=\int_{G(2,2n+2)} \ e(G(2,2n+2))=2n+2.$$
From (1) and
$$p(G(2, 2n+2)) =\prod_{\alpha=2}^{n+1} \
  (1+ 2(x_1^2+y_\alpha^2)+ (x_1^2-y_\alpha^2)^2),$$
we can prove (2).

Restricting the Euler class $e(E)$ on $G(2, k+2)$, we have
$$i^*e(E) =\frac 1{2\pi}\sum\limits_{\alpha=3}^{k+2} \ \omega_1^\alpha\wedge\omega_2^\alpha.$$
Then $$i^*e^{k}(E) =\frac {k!}{(2\pi)^{k}}
 \omega_1^3\wedge\omega_2^3\wedge\cdots\wedge \omega_1^{k+2}\wedge\omega_2^{k+2}=
 \frac {k!}{(2\pi)^{k}}\mbox{d}\, V_{G(2,k+2)},$$
where $\mbox{d}\, V_{G(2,k+2)}$ is the volume element of
$G(2,k+2)$. Then
$$\int_{G(2,k+2)} \ e^k(E(2,2n+2)) =2, \ \ \ k=1,\cdots,2n.$$

Restricting on $G(1, 2n+1)$, $\omega_2^\alpha=0$, we have
$$e(F(1,2n+1))=j^*e(F(2, 2n+2))=\frac {(2n)!}{2^{2n}n!\pi^{n}}
\omega_{1}^3\wedge\omega_{1}^4\wedge\cdots \wedge\omega_{1}^{2n+2}. $$
It is easy to see that $e(F(1,2n+1))$ is the Euler class of the tangent bundle of
$S^{2n}=G(1, 2n+1)$ and
$\omega_{1}^3\wedge\omega_{1}^4\wedge\cdots \wedge\omega_{1}^{2n+2}$ is
the volume element. Hence by Gauss-Bonnet formula or by direct computation, we have
$$\int_{G(1, 2n+1)} \ j^*e(F(2, 2n+2))= 2.$$

Furthermore, from $\omega_2^\alpha|_{G(1, 2n+1)}=0$ and $\Omega_{\alpha\beta}|_{G(2, n+2)}=0$
for $\alpha,\beta >n+2$,
 we have
$$\int_{G(1, 2n+1)} \ j^*e^n(E(2, 2n+2))=0, \ \ \ \int_{G(2, n+2)} \
i^*e(F(2, 2n+2))= 0.$$
By $p_1^k(E(2, 2n+2))=e^{2k}(E(2, 2n+2))$, we have
$$\int_{G(2,2k+2)} \
p_1^k(E(2, 2n+2))=2. \ \ \ \ \Box$$

The Poincar\' {e}  polynomial of $ G(2, 2n+2)$ is
$$p_t(G(2, 2n+2))=1+t^2+\cdots +t^{2n-2} + 2 t^{2n}+t^{2n+2} +\cdots + t^{4n}.$$
By Theorem 5.2, we have

(1) For $k\not = n$,  $e^k(E)\in H^{2k}(G(2,2n+2)), \ G(2,k+2)\in H_{2k}(G(2,2n+2))$ are the generators
 respectively;

(2) $e^n(E), e(F)\in H^{2n}(G(2,2n+2))$ and $G(2,n+2),
G(1,2n+1)\in H_{2n}(G(2,2n+2))$ are the generators.

The characteristic classes $e^k(E), e(F)$ and the submanifolds
$G(2,k+2), G(1,2n+1)$ are integral cohomology and homology classes respectively.
But they need not be the generators of the integral  cohomology and homology groups. For example,
when $k\not= n$, from
$\int_{G(2,k+2)} \ e^k(E)=2$ we know that $[G(2,k+2)]\in H_{2k}(G(2,2n+2), {\bf Z}), \
e^k(E)\in H^{2k}(G(2,2n+2), {\bf Z})$ can not be generators simultaneously.
Now we compute $\int_{{\bf C} P^k} \ i^*e^k(E(2, 2n+2))$ and
$\int_{{\bf C} P^n} \ i^*e(F(2, 2n+2))$.

Let $J$ be a complex structure on ${{\bf R}}^{2n+2}$,
${{\bf R}}^{2k+2}\subset {{\bf R}}^{2n+2}$ be invariant under the action of $J$ and
${\bf C} P^k= \{e_1Je_1 \ | \ e_1\in S^{2k+1}\}$. Let
$e_1,e_2=Je_1,e_{2\alpha-1}, e_{2\alpha}=Je_{2\alpha-1}, \ \alpha=2,3, \cdots,k+1$, be local Hermitian
orthonormal frame fields on ${{\bf R}}^{2k+2}$.
By $\mbox{d}e_2=J\mbox{d}e_1$ we have
$\omega_1^{2\alpha -1}=\omega_2^{2\alpha}, \ \ \omega_1^{2\alpha }=-\omega_2^{2\alpha-1},$ then
$$\mbox{d} (e_1\wedge e_2) = \sum\limits_{\alpha=2}^{k+1} \ \omega_1^{2 \, \alpha -1}
(E_{1 \,  2\alpha-1} + E_{2 \,  2\alpha})+
\sum\limits_{\alpha=2}^{k+1} \ \omega_1^{2\alpha }(E_{1  \, 2\alpha} - E_{2 \,  2\alpha-1}). $$
The oriented volume element of ${\bf C} P^k$ is
$\mbox{d}  V= 2^k \omega_1^3\wedge \omega_1^4\wedge\cdots \wedge\omega_1^{2k+2}.$

Let $i\colon\; {\bf C} P^k\to G(2, 2n+2)$ be inclusion,
we have
$$i^*e^k(E) = (-1)^k
\frac {k!}{\pi^k}\, \omega_1^3\wedge \omega_1^4\wedge\cdots \wedge\omega_1^{2k+2}.$$
By Proposition 2.3,
$$\int_{{\bf C} P^k} \ i^*e^k(E)=(-1)^k.$$

 $J$ induces a complex structure on
the induced bundle $i^*F \to {\bf C}P^n$. Let $F_{\bf C}$
be the complex vector bundle formed by the $(1,0)$-vectors of $i^*F\otimes  {\bf C}$.
By $i^*e(F)= c_n(F_{\bf C})$  (see [20]), we can show
$$\int_{{\bf C}P^n} \ i^*e(F)=\int_{{\bf C}P^n} \ c_n(F_{\bf C}) =1,$$
see also Chern [2].

Let $\bar J$ be a complex structure on ${{\bf R}}^{2k+2}$, the orientation given by
$\bar J$ is opposite to that of $J$.
 Let $\overline {{\bf C} P}^k=\{v\wedge \bar Jv \ | \ v\in S^{2k+1}\}$
be the complex projective space.
The orientation on the  vector bundle $E(2, 2n+2)|_{\overline {{\bf C} P}^k}$
are given by $v, \bar Jv$, we have
$$\int_{\overline {{\bf C} P}^k} \ i^*e^k(E)=(-1)^k.$$
By $p_k(E(2, 2n+2))= e^{2k}(E(2, 2n+2))$, we have
$$\int_{{{\bf C} P}^{2k}} \ i^*p_k(E(2, 2n+2))=1.$$

Let $\widetilde F_{\bf C}$
be the complex vector bundle formed by the $(1,0)$-vectors of
$F\otimes  {\bf C}|_{\overline {{\bf C} P}^n}$. The
orientation on realization vector bundle of $\widetilde  F_{\bf C}$ given by $\bar J$ is
opposite to  that of
$F|_{\overline {{\bf C} P}^n}$. Hence $e(F|_{\overline {{\bf C} P}^n})=
 -c_n(\widetilde  F_{\bf C})$ and we have
$$\int_{\overline {{\bf C} P}^n} \ e(F)=-\int_{\overline{{\bf C}P}^n} \ c_n(\widetilde F_{\bf C}) =-1.$$
These prove

{\bf Proposition 5.3} \ (1) When $k < n$, we have
$$[G(2,k+2)]=2(-1)^k[{\bf C}P^k]\in H_{2k}(G(2,2n+2));$$

(2) In the homology group $ H_{2n}(G(2,2n+2)) $, we have
$$[G(2,n+2)]=(-1)^n([{\bf C}P^n] +[\overline{{\bf C}P}^n]),$$
$$[G(1,2n+1)]=[{\bf C}P^n] -[\overline{{\bf C}P}^n].$$

For Grassmann manifold $G(2,2n+3)$, by the splitting principle of the characteristic classes,
 we can assume that there are oriented orthonormal
sections $s_{1},s_{2}$
and $t_{3},t_{4},\cdots,t_{2n+2}, t_{2n+3}$ of vector bundle $E(2k, 2n+3)$ and  $F(2, 2n+3)$ respectively,
such that
 $$\frac 1{2\pi}\nabla^2
\left (
  \begin{array}{ccccc}
      s_{1}\\
      s_{2}\\
        \end{array}
  \right )
  =
   \left (
  \begin{array}{ccccc}
    0 & -x\\
   x  & 0\\
     \end{array}
  \right )
  \left (
  \begin{array}{ccccc}
      s_{1}\\
      s_{2}\\
         \end{array}
  \right ), $$
$$\frac 1{2\pi}\nabla^2
\left (
  \begin{array}{cccccc}
      t_{3}\\
      t_{4}\\
      \vdots\\
      t_{2n+1} \\
      t_{2n+2} \\
      t_{2n+3}
  \end{array}
  \right )
  =
   \left (
  \begin{array}{cccccccc}
    0 & -y_{2} & & & \\
   y_{2}  & 0 & & \\
   & & \ddots & & & \\
  & & &  0 & -y_{n+1} &\\
  & & &  y_{n+1}  & 0 & \\
  & & &  & & & 0 \\
  \end{array}
  \right )
  \left (
  \begin{array}{ccccc}
      t_{3}\\
      t_{4}\\
      \vdots\\
      t_{2n+1} \\
      t_{2n+2} \\
      t_{2n+3}
  \end{array}
  \right ). $$
 The total Pontrjagin classes of $F$ are
 $p(F)=  \prod\limits_{\alpha=2}^{n+1} \ (1+ y_\alpha^2).$

  $ s_{1}\otimes  t_{2\alpha-1},
   s_{2}\otimes  t_{2\alpha-1},
   s_{1}\otimes  t_{2\alpha},
   s_{2}\otimes  t_{2\alpha}$ and $s_{1}\otimes  t_{2n+3},
   s_{2}\otimes  t_{2n+3}$ are orthonormal sections of
$E\otimes F\cong TG(2k,2n)$, they also give an orientation of $E\otimes F$.
The curvature of $E\otimes F$ are defined by
 $$\frac 1{2\pi}\nabla^2
 \left (
 \begin{array}{ccccccccc}
   s_{1}\otimes  t_{2\alpha-1}\\
   s_{2}\otimes  t_{2\alpha-1}\\
   s_{1}\otimes  t_{2\alpha}\\
   s_{2}\otimes  t_{2\alpha}
 \end{array}
 \right )
 =  \left (
 \begin{array}{ccccccccc}
    0 & -x & -y_\alpha & 0\\
    x & 0 & 0 &  -y_\alpha \\
   y_\alpha & 0 & 0 & -x\\
   0 & y_\alpha & x & 0
 \end{array}
 \right )
 \left (
 \begin{array}{ccccccccc}
 s_{1}\otimes  t_{2\alpha-1}\\
   s_{2}\otimes  t_{2\alpha-1}\\
   s_{1}\otimes  t_{2\alpha}\\
   s_{2}\otimes  t_{2\alpha}
 \end{array}
 \right ), $$
 $$\frac 1{2\pi}\nabla^2
 \left (
 \begin{array}{ccccccccc}
   s_{1}\otimes  t_{2n+3}\\
   s_{2}\otimes  t_{2n+3}
   \end{array}
 \right )
 =  \left (
 \begin{array}{ccccccccc}
    0 & -x \\
    x & 0 \\
  \end{array}
 \right )
 \left (
 \begin{array}{ccccccccc}
 s_{1}\otimes  t_{2n+3}\\
   s_{2}\otimes  t_{2n+3} \end{array}
 \right ). $$
Hence the Euler class of Grassmann manifold $G(2,2n+3)$ is
 $$e(TG(2, 2n+3)) =e(E\otimes F)= x \prod_{\alpha=2}^{n+1} \
  (x^2-y_\alpha^2)=(n+1)e^{2n+1}(E).$$

The odd dimensional homology groups of  Grassmann manifold $G(2,2n+3)$ are also trivial,
the even dimensional homology groups are one dimensional.
The Euler-Poincar\'{e} number is
$\chi(G(2,2n+3))=2n+2.$

Similar to the case of $G(2,2n+2)$, we have

{\bf  Theorem 5.4} \
(1)  The Pontrjagin classes of  $F(2,2n+3)$ and
 $TG(2,2n+3)$ can all be represented
by the Euler class  $e(E(2,2n+3))$;

 (2) $e(TG(2, 2n+3)) =(n+1)e^{2n+1}(E(2, 2n+3));$

 (3) $ \int_{G(2,k+2)} \ e^{k}(E(2, 2n+3))=2, \ \ k=1,\cdots,2n+1;$

(4) $ \int_{{{\bf C} P}^{k}} \ e^{k}(E(2, 2n+3))=(-1)^k, \ \ k=1,\cdots,n+1.$

As is well-known, the Chern, Pontrjagin and Euler classes are all integral cocycles.
Let $D\colon\; H^k(G(2, N),{{\bf Z}}) \to H_{2N-4-k}(G(2, N),{{\bf Z}})$ be the Poincar\' {e} duality.
The following theorem gives  the structure of the integral homology and coholomogy of $G(2,N)$.

{\bf   Theorem 5.5} \ (1) When $2k+2 < N$, $[{\bf C}P^k]$ and
 $e^k(E(2, N))$ are the generators of $H_{2k}(G(2,N), {{\bf Z}})$ and
 $H^{2k}(G(2,N), {{\bf Z}})$ respectively;

(2) When $2k+2 > N$, $[G(2,k+2)]$ and
$\frac 12e^k(E(2, N))$
are the generators of $ H_{2k}(G(2,N), {{\bf Z}})$ and
$H^{2k}(G(2,N),{{\bf Z}})$ respectively;

(3)   When $2k+2 < N$, $D(e^k(E(2, N)))=[G(2,N-k)]$;
when $2k+2 > N$,  $D(\frac 12 e^k(E(2, N)))=(-1)^{n-k}[{\bf C}P^{n-k-2}]$;

(4) $[{\bf C}P^n], [\overline{{\bf C}P}^n]$
and $\frac 12 \{(-1)^ne^n(E(2, 2n+2))\pm e(F(2, 2n+2))\}$
are  generators of $ H_{2n}(G(2,2n+2),{{\bf Z}})$ and $H^{2n}(G(2,2n+2), {{\bf Z}})$ respectively. Furthermore,
$$D(\frac 12 [(-1)^ne^n(E(2, 2n+2))+ e(F(2, 2n+2))]) =[{\bf C}P^n],$$
$$D(\frac 12 [(-1)^ne^n(E(2, 2n+2))- e(F(2, 2n+2))]) =[\overline{{\bf C}P}^n].$$

{\bf Proof} \ As is well-known, the Euler classes and Pontrjagin classes  are
harmonic forms and are integral cocyles, their product are also harmonic forms, see [7], [13].
When $2k+2 < N$, from
$\int_{{\bf C} P^k} \ e^k(E(2, N))=(-1)^k$ we know
${\bf C}P^k\in
H_{2k}(G(2,N), {{\bf Z}})$ and
 $e^k(E(2, N)) \in H^{2k}(G(2,N), {{\bf Z}})$ are  generators respectively.

By simple computation, we have
$$e^k(E(2,N))=\frac {k!}{(2\pi)^k}\sum\limits_{\alpha_1<\cdots < \alpha_k} \ \omega_1^{\alpha_1}
\omega_2^{\alpha_1}\cdots\omega_1^{\alpha_k}
\omega_2^{\alpha_k},$$
$$a=(e^k(E(2,N)),e^k(E(2,N)))=\frac {(k!)^2}{(2\pi)^{2k}}C_{N-2}^kV(G(2,N)),$$
\begin{eqnarray*}\frac 1a*e^k(E(2,N))& =  &
\frac {k!}{a(2\pi)^k}\sum\limits_{\beta_1<\cdots <\beta_{N-k-2}} \ \omega_1^{\beta_1}
\omega_2^{\beta_1}\cdots\omega_1^{\beta_{N-k-2}}
\omega_2^{\beta_{N-k-2}} \\
&=& \frac {(N-k-2)!}{2(2\pi)^{N-k-2}}\sum\limits_{\beta_1<\cdots <\beta_{N-k-2}} \ \omega_1^{\beta_1}
\omega_2^{\beta_1}\cdots\omega_1^{\beta_{N-k-2}}
\omega_2^{\beta_{N-k-2}} \\
&=&\frac 12 e^{N-k-2}(E(2,N)).\end{eqnarray*}
By  Theorem 3.1, $\frac 12 e^{N-k-2}(E(2,N))$ is a generator of $ H^{2N-2k-4}(G(2,N), {{\bf Z}})$.
By $\int_{G(2,N-k)} \ \frac 12 e^{N-k-2}(E(2,N))=1$ we know that
$G(2,N-k)\in  H_{2N-2k-4}(G(2,N), {{\bf Z}})$ is a generator and $D(e^k(E(2, N)))=G(2,N-k)$.
These proves (1), (2), (3) of the Theorem.

Let $[S_1],[S_2]$ be generators of $H_{2n}(G(2,2n+2),{{\bf Z}})$ and harmonic forms
$\xi_1,\xi_2$ be generators of $ H_{2n}(G(2,2n+2),{{\bf Z}})$, they satisfy
$\int_{S_i} \ \xi_j =\delta_{ij}.$
There are integers $a_{ij}, n_{ij}$ such that
$$\left(\begin{array} {cc}  e^n(E(2,2n+2)) \\ e(F(2,2n+2)) \end{array} \right)
= \left(\begin{array} {cc} n_{11} & n_{12} \\
n_{21} & n_{22} \end{array} \right)\left( \begin{array} {c} \xi_1 \\  \xi_2
\end{array} \right), $$
$$ ({\bf C}P^n, \overline{{\bf C}P}^n)
= (S_1,  S_2 )\left(\begin{array} {cc} a_{11} & a_{12} \\
a_{21} & a_{22} \end{array} \right). $$
Then
$$\left(\begin{array} {cc} n_{11} & n_{12} \\
n_{21} & n_{22} \end{array} \right)\left(\begin{array} {cc} a_{11} & a_{12} \\
a_{21} & a_{22} \end{array} \right)=\left(\begin{array} {cc} (-1)^n & (-1)^n \\
1 & -1 \end{array} \right),$$
we have $\det (a_{ij})=\pm 1$ or $\det (n_{ij})=\pm 1$.

Assuming $\det (n_{ij})=\pm 1$, then $e^n(E(2,2n+2)), e(F(2,2n+2))$ are  the generators of
$ H^{2n}(G(2, 2n+2), {{\bf Z}})$, we can assume $\xi_1=e^n(E(2,2n+2)),\xi_2=e(F(2,2n+2))$.
It is easy to see $$(e^n(E(2,2n+2)),e(F(2,2n+2)))=0,$$ $$(e^n(E(2,2n+2)),e^n(E(2,2n+2)))=
(e(F(2,2n+2)),e(F(2,2n+2)))=2.$$ By Theorem 3.1,
$\frac 12 e^n(E(2,2n+2)), \ \frac 12 e(F(2,2n+2))$ are also  the generators of
$ H^{2n}(G(2,2n+2),{{\bf Z}})$. This contradict to
 $e^n(E(2,2n+2)), e(F(2,2n+2))$ be  the generators.

Then we must have $\det (n_{ij})=\pm 2$ and $\det (a_{ij})=\pm 1$. These shows
$ {\bf C}P^n, \overline{{\bf C}P}^n$ are
generators of $ H_{2n}(G(2,2n+2),{{\bf Z}})$,
$\frac 12 \{(-1)^ne^n(E(2, 2n+2))+ e(F(2, 2n+2))\}$ and
$\frac 12 \{(-1)^ne^n(E(2, 2n+2))- e(F(2, 2n+2))\}$
are generators of  $ H^{2n}(G(2,2n+2),{{\bf Z}})$.
The Poincar\' {e} duals of these  generators are easy to compute.  \ \ \ \ $\Box$

Finally we give some applications to conclude   this section.

Let $f\colon\; M  \to  {{\bf R}}^N$ be an immersion of an oriented compact surface $M$,
$g\colon\; M \to G(2, N)$ the induced Gauss map, $g(p)=T_pM$. Then
$e(M)=g^*(e(E(2, N))$ is the Euler class of $M$.
Let $[M]\in H_2(M)$ be the fundamental class of $M$. When $N\not= 4$, we have
$$g_*[M]= \frac 12\chi(M)[G(2,3)]= \chi(M)[-{\bf C}P^1]\in H_2(G(2,N)).$$

In [15], [16], we have shown there is a fibre bundle $\tau\colon\;   G(2, 8) =G(6,8)
\to S^6$ with fibres ${\bf C}P^3$, where $S^6=\{ v\in S^7 \ | \ v\perp \bar e_1=(1,0,\cdots,0)\}$,
 $\tau^{-1}(\bar e_2)=
\{v\wedge Jv \ | \ v\in S^7\}$, $\bar e_2=(0,1,0,\cdots,0)$.
On the other hand, the map $f(v)=\bar e_1\wedge v$ gives a section of $\tau$.
Let $\mbox{d}V$ be the volume form on $ S^6$ such that $\int_{S^6} \ \mbox{d}V=1$.
It is easy to see
$$[\tau^*\mbox{d}V] =\frac 12 e^3(E(2,8)) + \frac 12 e(F(2,8)).$$

Let  $\varphi\colon\; M\to {{\bf R}}^8$
be an immersion of an oriented compact $6$ dimensional manifold,
$g\colon\; M\to  G(6, 8) =G(2,8)$ be Gauss map. Then
$e(M)=g^*e(F(2,8))$ is the Euler class of tangent bundle of $M$,
$e(T^\perp M)=g^*e(E(2,8))$ is the Euler class of normal bundle of $M$.
$$\int_M \ (\tau\circ g)^*\mbox{d}V= \int_{M} \ g^*[\frac 12 e^3(E(2,8)) + \frac 12 e(F(2,8))]
= \frac 12\int_{M} \  e^3(T^\perp M) + \frac 12 \chi (M) $$
is the degree of the map $\tau\circ g\colon\; M\to S^6$.
If $\varphi$ is an imbedding,  $e(T^\perp M)=0$,
see Milnor, Stasheff [12], p.120.

Let $J, \ \bar J$ be two complex structures on ${{\bf R}}^4$, with
orthonormal basis $e_1,e_2, e_3,e_4$,
$$Je_1=e_2, \ Je_2=-e_1, \ Je_3=e_4, \ Je_4=- e_3;$$
$$\bar Je_1=-e_2, \ \bar Je_2=e_1, \
\bar Je_3=e_4,\ \bar Je_4=-e_3.$$
For any unit vector $v=\sum \ v_ie_i$, we have
$$ vJv+ * vJv  = e_1e_2 + e_3e_4,$$
\begin{eqnarray*} vJv- * vJv & = & (v_1^2 + v_2^2-v_3^2 - v_4^2)(e_1e_2-e_3e_4) \\
 & & + 2(v_1v_3 + v_2v_4)(e_1e_4
-e_2e_3) + 2(v_2v_3-v_1v_4 )(e_1e_3
+e_2e_4); \end{eqnarray*}
$$ v\bar Jv- * v\bar Jv  = - e_1e_2 + e_3e_4,$$
\begin{eqnarray*} v\bar Jv + * v\bar Jv & = & (-v_1^2 - v_2^2+v_3^2 + v_4^2)(e_1e_2+e_3e_4) \\
 & & + 2(v_1v_3 - v_2v_4)(e_1e_4
+e_2e_3) - 2(v_1v_4+v_2v_3 )(e_1e_3
-e_2e_4).\end{eqnarray*}
These shows ${\bf C}P^1$, $\overline{{\bf C}P}^1$ are two spheres in
$G(2,4)\approx S^2(\frac {\sqrt 2}2)\times
S^2(\frac {\sqrt 2}2)$ where the decomposition is given by Hodge star operator.
Let $f\colon\; M  \to  {{\bf R}}^4$ be an imbedding of an oriented surface,
$g\colon\; M \to G(2, 4)$ the Gauss map. Then we have
$$g_*[M]= \frac 12\chi(M)[G(2,3)]=-\frac 12\chi(M)[{\bf C}P^1]- \frac 12\chi(M)
[\overline{{\bf C}P}^1].$$
See the work of Chern and Spanier [1].

\vskip 1cm \noindent{\bf \S 6. \  The case of $G(3,6)$ } \vskip 0.3cm

The Poincar\' {e}  polynomial of $ G(3, 6)$ is
$p_t(G(3, 6))=1+t^4+ t^{5} + t^{9}.$
To study the real homology and coholomogy of
 Grassmann manifold $G(3,6)$ we need only consider the dimension $4$ and $5$.

Let $i\colon\; G(2,4)\to G(3,6)$ be inclusion defined naturally. It is easy to see that $i^*
p_1(E(3,6))=p_1(E(2,4))=e^2(E(2,4))$, then
$$\int_{G(2,4)} \ p_1(E(3,6))=2.$$

As \S 2, let $SLAG = \{G(\bar e_1\bar e_3\bar e_5) \ | \ G\in SU(3)\subset SO(6)\}$ be a subspace of $G(3,6)$,
$e_i=G(\bar e_{i}), e_{i+1}=G(\bar e_{i+1})=Je_{i}$ be $SU(3)$-frame fields, \ $i=1,3,5$.
Restricting the coframes $\omega_A^B=\langle \mbox de_A, e_B\rangle $ on $SLAG$ we have
$$\omega_i^j=\omega_{i+1}^{j+1}, \ \omega_i^{j+1}= -\omega_{i+1}^{j}, \ i,j=1,3,5 \ \ \mbox{and} \ \
\omega_1^{2}+\omega_3^{4} + \omega_5^{6}=0.$$
By the proof of Proposition 2.4, we have $\mbox {d}V_{SLAG} = 2^{\frac 32}\sqrt 3\,
\omega_1^4\omega_1^6\omega_3^6\omega_1^2
\omega_3^4$ and
$V(SLAG) \\ = \sqrt {\frac 32} \, \pi^3$. Let $G(3,6)$ be generated by $e_1e_3e_5$ locally,
the first  Pontrjagin class of canonical vector bundle
$E(3,6)$ is
$$p_1(E(3,6)) = \frac 1{4\pi^2} [(\Omega_{13})^2 + (\Omega_{15})^2 +  (\Omega_{35})^2 ],$$
where $\Omega_{ij}=-\sum\limits_{\alpha} \  \omega_{i}^\alpha \wedge\omega_{j}^\alpha, \ \alpha=2,4,6$.
By computation we have
$$*p_1(E(3,6))|_{SLAG}=\frac {\sqrt 6}{4\pi^2}\mbox{d} V_{SLAG},$$
$$a=(p_1(E(3,6)),p_1(E(3,6)))=\frac {3\cdot 4\cdot 3}{(2\pi)^4} V(G(3,6))=\frac 32\pi.$$

From $\int_{G(2,4)} \ p_1(E(3,6))=2$, we know that
$p_1(E(3,6))$ or $\frac 12 p_1(E(3,6))$ is a generator
of $ H^{4}(G(3,6),{{\bf Z}})$. If $p_1(E(3,6))$ is a generator, by Theorem 3.1,
$\frac 1a*p_1(E(3,6))$ is a generator of $ H^{5}(G(3,6),{{\bf Z}})$, but
$$\int_{SLAG} \  \frac 1a*p_1(E(3,6))=\int_{SLAG} \ \frac 1{\sqrt 6 \pi^3}\mbox {d} V_{SLAG} =\frac 12.$$
Then $\frac 12 p_1(E(3,6))$ is a generator
of $ H^{4}(G(3,6),{{\bf Z}})$ and  $\int_{SLAG} \  \frac 4a*\frac 12 p_1(E(3,6))=1$.

We have proved following theorem

{\bf  Theorem 6.1} \ (1) $\frac 12 p_1(E(3,6))\in  H^{4}(G(3,6),{{\bf Z}})$ is a generator and
its Poincar\' {e} dual  $[SLAG]$ is a generator of $H_5(G(3,6),{{\bf Z}})$;

(2) $\frac 4{3\pi} *p_1(E(3,6))\in  H^{5}(G(3,6),{{\bf Z}})$ is a generator and
its Poincar\' {e} dual  $[G(2,4)]$ is a generator of $H_4(G(3,6),{{\bf Z}})$.

Let $\bar e_1,\cdots,\bar e_6$ be a fixed orthonormal basis of ${\bf R}^6$, $G\in SO(3)$  acts
on the subspace generated by $\bar e_4,\bar e_5, \bar e_6$, denote $e_4=G(\bar e_4),
e_5=G(\bar e_5), e_6=G(\bar e_6)$.
As [7], let $PONT$ be the set of elements
$$ (\cos t\bar e_1 +\sin t e_4) (\cos t\bar e_2 +\sin t e_5)
 (\cos t\bar e_3 +\sin t e_6), \ \ t\in [0,\frac {\pi}2].$$
$PONT$ is a calibrated submanifold (except two points correspond to $t=0,\frac {\pi}2$) of the first
Pontrjagin form $p_1(G(3,6))$. By moving frame we can show
$$\int_{PONT} \ p_1(G(3,6))=2 \ \ \mbox{and} \ \ V(PONT)=\sqrt {\frac 23}\, V(G(2,4))=\frac {4\sqrt 6}3\pi^2.$$
Then $4$-cycle $PONT$ is homologous to the $4$-cycle $G(2,4)$ inside $G(3,6)$.

\vskip 0.5cm \noindent{\bf \S 7.  The case of $G(3,7)$ } \vskip 0.3cm

The Poincar\' {e}  polynomial of $ G(3, 7)$ is
$p_t(G(3, 7))=1+ 2t^4+ 2t^{8} + t^{12}.$

Let  $\bar e_1, \bar e_2, \cdots,   \bar e_8$ be a fixed orthonormal basis of ${\bf R}^8$ and
 ${\bf R}^7$ be a subspace generated by $\bar e_2, \cdots,   \bar e_8$.
The oriented Grassmann manifold $G(3,7)$ is the set of subspaces of ${\bf R}^7$.

Let $E=E(3,7)$ and $F=F(3,7)$ be canonical vector bundles on Grassmann manifold $G(3,7)$.
As \S 4, we can show
$$p_1(F)=-p_1(E), \ p_2(F)=e^2(F)=p_1^2(E), \ e(E\otimes F)=e(TG(3,7))=3 e^3(F).$$
By $\int_{G(3,7)} \ e(TG(3,7))=\chi(G(3,7))=6$ we have
$$\int_{G(3,7)} \ e^3(F)=2.$$
By inclusion $G(2,6)\subset G(3,7)$,  ${\bf C}P^2$ and $\overline {{\bf C}P}^2$
can be imbedded in $G(3,7)$.

{\bf  Lemma 7.1} \ $\int_{{\bf C}P^2} \ p_1(E) = \int_{\overline {{\bf C}P}^2} \ p_1(E)=1,
 \ \ \int_{{\bf C}P^2} \ e(F)=- \int_{\overline {{\bf C}P}^2} \ e(F)=1.$

{\bf  Proof} \ By $p_1(E(3,7))|_{{\bf C}P^2}=p_1(E(2,6))|_{{\bf C}P^2}=
e^2(E(2,6))|_{{\bf C}P^2}$ and the results of \S 5, we have
$\int_{{\bf C}P^2} \ p_1(E)=1$. The other equalities can be proved similarly.
 \ \ \ \ $\Box$

Then
$$\int_{{\bf C}P^2} \ \frac 12(p_1(E)+e(F))
= \int_{\overline {{\bf C}P}^2} \ \frac 12(p_1(E)-e(F))=1,$$
$$\int_{{\bf C}P^2} \ \frac 12(p_1(E)-e(F))
= \int_{\overline {{\bf C}P}^2} \ \frac 12(p_1(E)+e(F))=0,$$
hence ${\bf C}P^2,  \overline {{\bf C}P}^2\in H_4(G(3,7) )$
and $\frac 12(p_1(E)+e(F)), \
\frac 12(p_1(E)-e(F))\in H^4(G(3,7) )$ are generators.

Let $e_2,e_3,e_4,\cdots,e_{8}$ be oriented orthonormal frame fields on ${{\bf R}}^{7}$,
$G(3,7)$ be generated by $e_2\wedge e_3\wedge e_4$ locally. Euler class of $F$ and first
Pontrjagin class of $E$ can be represented by
\begin{eqnarray*} e(F) &=&\frac {1}{2(4\pi)^2}\sum\limits \ \varepsilon(\alpha_1\alpha_2\alpha_{4}\alpha_{4})
\Omega_{\alpha_1\alpha_2}\wedge \Omega_{\alpha_{3}\alpha_{4}} \\
&=& \frac 1{4\pi^2} \sum\limits_{i,j=2}^4 \ (\omega_5^i\omega_6^i\omega_7^j\omega_8^j
-  \omega_5^i\omega_7^i\omega_6^j\omega_8^j+\omega_5^i\omega_8^i\omega_6^j\omega_7^j ),\end{eqnarray*}
$$p_1(E) = \frac 1{4\pi^2} [(\Omega_{23})^2 + (\Omega_{24})^2 +  (\Omega_{34})^2 ].$$

{\bf  Lemma 7.2} \ (1) $* p_1(E)= \frac 45 \pi^2 p_1(E)e(F), \ *e(F)=\frac 12 \pi^2 e^2(F)$;

(2) $(p_1(E),e(F))=0, \ a= (p_1(E),p_1(E))=\frac 85\pi^2, \
b= (e(F),e(F))=\pi^2$.

{\bf Proof} \  $* p_1(E), \ p_1(E)e(F)$ and $*e(F), \ e^2(F)$
 are the harmonic forms on $G(3,7)$. $* p_1(E)= \frac 45 \pi^2 p_1(E)e(F)$ follows from the equalities
 such as
\begin{eqnarray*} *\omega_2^5\omega_3^5\omega_2^6\omega_3^6
&= & \omega_2^7\omega_3^7\omega_2^8\omega_3^8\omega^4_5\omega^4_6\omega^4_7\omega^4_8 \\
&= & \omega_3^7\omega_4^7\omega_3^8\omega_4^8\omega^4_5\omega^4_6\omega^2_7\omega^2_8 \\
&= & \omega_2^7\omega_4^7\omega_2^8\omega_4^8\omega^4_5\omega^4_6\omega^3_7\omega^3_8. \end{eqnarray*}
The proof of (2) is a direct computation. \ \ \ \ $\Box$

To study $G(3,7), G(3,8)$ and $ G(4,8)$, we shall use Clifford algebras.

 Let $C\ell_8$ be the Clifford algebra associated to
 the Euclidean space ${\bf R}^8$. Let
  $\bar e_1,   \bar e_2, \cdots,   \bar e_8$ be a fixed orthonormal basis of ${\bf R}^8$,
the Clifford product be determined by the relations: $\bar
e_B\cdot \bar e_C+\bar e_C\cdot \bar e_B =-2\delta_{BC}, \ \ B,  C=1,  2,
\cdots,  8$.  Define the subspace  $ V=V^+\oplus V^-$ of $C\ell_8$ by
$  V^+=C\ell^{even}_8\cdot A, V^-=C\ell^{odd}_8\cdot A$,
where
 $$A=\frac 1{16} \mbox {Re} \, [(\bar
e_1 + \sqrt{-1}\bar e_2)\cdots(\bar e_7 + \sqrt{-1}\bar e_8) (1 + \bar
e_1\bar e_3\bar e_5\bar e_7)].  $$
The space $V=V^+\oplus V^-$ is an
irreducible module over $C\ell_8$. The spaces $V^+$ and $V^-$ are
generated by $\bar e_1\bar e_BA$ and $\bar e_B A$ respectively,
$B=1, \cdots, 8$, see [14], [15].

Let $Spin_7=\{ G \in SO(8) \ | \ G(A)=A \}$ be the isotropy group of
$SO(8)$ acting on $A$. The group $Spin_7$ acts
on $G(2,8), G(3,8)$ and $S^7$ transitively. $G_2=\{ G\in Spin_7 \ | \ G(\bar e_1)=\bar e_1 \}$
is a subgroup of $Spin_7$.

 The Grassmann manifold $G(k, 8)$ can  be viewed as  subset
of Clifford algebra $C\ell_8$ naturally. Then for any $\pi\in G(k,8)$, there
is $v\in{\bf R}^8$ such that $\pi A=\bar e_1vA  $ or $\pi A=vA $ according to
the number $k$ being even or odd, $|v|=1$. Thus we have  maps
$G(k, 8)\to  S^7, \ \pi \mapsto v$.
Since $Spin_7$ acts on $G(3,8)$ transitively, from $\bar
e_2\bar e_3\bar e_4 A=\bar e_1 A$ we have $G(\bar
e_2\bar e_3\bar e_4) A=G(\bar e_1) A$ for any $G\in Spin_7$. This shows the map
$\tau\colon\;  G(3, 8)\to  S^7, \ \tau (\pi)=v,$ is a fibre bundle and $v\perp \pi$,  see [15], [16].
Let $$ASSOC=\tau^{-1}(\bar e_1) = \{ \pi \in G(3,8) \ | \ \tau(\pi)=\bar e_1\}$$ be the fibre over $\bar e_1$.
The group $G_2$ acts on  $ASSOC$ transitively, we have
$ASSOC=\{ G(\bar e_2\bar e_3\bar e_4) \ | \ G\in G_2\}$. We can show the isotropy group
 $\{ G(\bar e_2\bar e_3\bar e_4)=\bar e_2\bar e_3\bar e_4 \ | \ G\in G_2\}$ is isomorphic to
the group $SO(4)$, then $ASSOC\approx G_2/SO(4)$.

Change the orientation of ${\bf R}^7$, let $\tilde A=\frac 1{16}\mbox {Re} \, [(\bar
e_1-\sqrt{-1}\bar e_2)(\bar e_3+\sqrt{-1}\bar e_4)\cdots(\bar e_7+\sqrt{-1}\bar e_8) (1+\bar
e_1\bar e_3\bar e_5\bar e_7)]$.
Define submanifold $\widetilde {ASSOC}= \{\pi\in G(3,8) \ | \ \pi \tilde A =\bar e_1 \}$ which
diffeomorphic to $ASSOC$.

{\bf  Lemma 7.3} \ $V(ASSOC)= \frac 65\pi^4$.

{\bf  Proof} \  Let $ \tilde e_1,\tilde e_2,
\cdots, \tilde e_8$ be $Spin_7$ frame fields on ${\bf R}^8$, the 1-forms
$\omega_B^C =$  $\langle \mbox {d}\tilde e_B, \tilde e_C\rangle $ satisfy (for proof, see [16] )
\begin{eqnarray*}
&\omega_1^2+\omega_3^4+\omega_5^6+\omega_7^8=0, \ &
\omega_1^3-\omega_2^4+\omega_6^8-\omega_5^7=0, \\ &
\omega_1^4+\omega_2^3+\omega_5^8+\omega_6^7=0,  \
&\omega_1^5-\omega_2^6+\omega_3^7-\omega_4^8=0, \\
&     \omega_1^6+\omega_2^5-\omega_3^8-\omega_4^7=0,  \
&\omega_1^7-\omega_2^8-\omega_3^5+\omega_4^6=0, \\
&  \omega_1^8+\omega_2^7+\omega_3^6+\omega_4^5=0.  \ & \quad
\end{eqnarray*}
Since $Spin_7$ acts on $G(3,8)$  transitively, $G(3,8)$ is locally generated by $\tilde e_2\tilde e_3
 \tilde e_4$. The volume element of $G(3,8)$ is $\mbox{d} V_{G(3,8)}=
\omega_2^1\omega_3^1\omega_4^1\omega_2^5\omega_3^5\omega_4^5\cdots \omega_2^8\omega_3^8\omega_4^8$.

Note that $A$ can be represented by $ Spin_7$ frames, that is
$$A=\frac 1{16}\mbox {Re} \, [(\tilde
e_1+\sqrt{-1}\tilde  e_2)\cdots(\tilde  e_7+\sqrt{-1}\tilde e_8) (1+\tilde
e_1\tilde  e_3\tilde  e_5\tilde e_7)].  $$

Let $ \tilde e_1=\bar e_1$ be a fixed vector, $ \bar e_1,\tilde e_2,
\cdots, \tilde e_8$ be $G_2$ frame fields on ${\bf R}^8$,
$ASSOC$ is locally generated by $\tilde e_2\tilde e_3\tilde e_4$ and
\begin{eqnarray*} \mbox {d}(\tilde e_2\tilde e_3\tilde e_4) &=& \sum\limits_{i=2}^3
\sum\limits_{\alpha=5}^8 \ \omega_i^\alpha E_{i\alpha} \\
& = & \omega_2^5 (E_{25} + E_{47}) +\omega_2^6 (E_{26} - E_{48})+\omega_2^7 (E_{27} - E_{45})
+\omega_2^8 (E_{28} + E_{46}) \\
& & +
\omega_3^5 (E_{35} + E_{46}) +\omega_3^6 (E_{36} - E_{45})+\omega_3^7 (E_{37} + E_{48})
+\omega_3^8 (E_{38} - E_{47}). \end{eqnarray*}
The metric on $ASSOC$ is
\begin{eqnarray*} \mbox {d}s^2
& = & 2(\omega^5_2)^2+ 2(\omega^8_3)^2-2\omega^5_2\omega^8_3
+  2(\omega^6_2)^2+ 2(\omega^7_3)^2-2\omega^6_2\omega^7_3 \\
&& + 2(\omega^7_2)^2+ 2(\omega^6_3)^2+2\omega^7_2\omega^6_3
+  2(\omega^8_2)^2+ 2(\omega^5_3)^2+2\omega^8_2\omega^5_3,\end{eqnarray*} with the volume form
 $$\mbox {d}V_{ASSOC}=9 \ \omega^5_2\omega^5_3\ \cdots\ \omega^8_2\omega^8_3.$$
The normal space of $ASSOC$ in $G(3,8)$ at $\tilde e_2\tilde e_3\tilde e_4$ is generated by
$$E_{21}, E_{31}, E_{41}, E_{25} - E_{47} -E_{38}, E_{26} + E_{48} -E_{37},
E_{27} + E_{45} +E_{36}, E_{28} - E_{46} +E_{35}.$$

The sphere $S^7$ is generated by $\tilde e_1$,  $\mbox{d} V_{S^7} =
\omega_2^1\omega_3^1\cdots\omega_1^8$ is the volume form. From
$$(E_{27}+E_{45}+E_{36})A= -3 \tilde e_8A, \ E_{21}A=-\tilde e_2A, \ \cdots,$$
we can compute the tangent map of $\tau\colon\;  G(3, 8)\to  S^7$,
$$\tau_*(E_{27}+E_{45}+E_{36})= -3 \tilde e_8, \ \tau_*(E_{21})=-\tilde e_2, \ \cdots.$$
Then we can compute the cotangent map $\tau^*$ and we have
$$\tau^*\mbox{d} V_{S^7}= \omega^1_2\omega_3^1\omega_4^1
(\omega_2^6-\omega_3^7+\omega_4^8)(-\omega_2^5+\omega_3^8+\omega_4^7)
(\omega_2^8+\omega_3^5-\omega_4^6)(\omega_2^7+\omega_3^6+\omega_4^5),$$
and
$$\mbox{d} V_{G(3,8)}=-\frac 19 \tau^*\mbox{d} V_{S^7}\cdot\mbox {d}V_{\tau^{-1}(\tilde e_1)}.$$
From $V(G(3,8))=\frac 2{45}\pi^8, \ V(S^7)=\frac 1{3}\pi^4$, we have $V(ASSOC)= \frac 65\pi^4$.
 \ \ \ \ $\Box$

{\bf  Lemma 7.4} \ $\int_{ASSOC} \ p_1^2(E)=\int_{ASSOC} \ p_1(E)e(F)=1,$

\qquad $\int_{\widetilde {ASSOC}} \ p_1^2(E)=-\int_{\widetilde {ASSOC}} \ p_1(E)e(F)=1.$

\noindent Then  $ASSOC, \widetilde {ASSOC}$
and $p_1^2(E),  p_1(E)e(F)$ are generators of  $H_8(G(3,7) )$ and $
H^8(G(3,7) )$ respectively. Furthermore, we have
$$[G(2,6)]=[ASSOC]+ [\widetilde {ASSOC}].$$

{\bf  Proof} \ From
$$ *\omega_2^5\omega_3^5\omega_2^6\omega_3^6|_{ASSOC}
= \omega_2^7\omega_3^7\omega_2^8\omega_3^8\omega^4_5\omega^4_6\omega^4_7\omega^4_8|_{ASSOC}=\frac 19
\mbox{d}V_{ASSOC},$$
$$ *\omega_2^7\omega_3^7\omega_2^8\omega_3^8|_{ASSOC}
= *\omega_2^5\omega_3^5\omega_2^7\omega_3^7|_{ASSOC}=
*\omega_2^6\omega_3^6\omega_2^8\omega_3^8|_{ASSOC}  =\frac 19
\mbox{d}V_{ASSOC},$$
$$ *\omega_2^6\omega_3^6\omega_2^7\omega_3^7|_{ASSOC}
= \omega_2^5\omega_3^5\omega_2^8\omega_3^8\omega^4_5\omega^4_6\omega^4_7\omega^4_8|_{ASSOC}=0,$$
$$ *\omega_2^5\omega_3^5\omega_2^8\omega_3^8|_{ASSOC} =0, $$
we have $\sum\limits_{\alpha, \beta} \  *\omega_2^\alpha\omega_3^\alpha\omega_2^\beta\omega_3^\beta|_{ASSOC}
=8\cdot \frac 19\mbox{d}V_{ASSOC}$, $$* p_1(E)|_{ASSOC}= \frac 1{4\pi^2}\cdot 3\cdot 8\cdot \frac 19
\mbox{d}V_{ASSOC}.$$
Then by $* p_1(E)= \frac 45 \pi^2 p_1(E)e(F)$, we have
$$\int_{ASSOC} \  p_1(E)e(F) = \int_{ASSOC} \ \frac 5{4 \pi^2}*p_1(E)=1.$$

The proof of $\int_{ASSOC} \ p_1^2(E)=1$ is similar.
 \ \ \ \ $\Box$

Change the orientation of ${\bf R}^7$, we have Euclidean space
 $\widetilde {{\bf R}}^7$. Let $\widetilde E, \widetilde F\to G(3,7)$ be canonical vector
bundles with respect to $\widetilde {{\bf R}}^7$. It is easy to see that $\widetilde E=E$,
but the orientations of $\widetilde F$ and $F$ are different. These shows
$$\int_{\widetilde {ASSOC}} \ p_1^2(E)=1, \ \
 \int_{\widetilde {ASSOC}} \ p_1(E)e(F)=-1.$$

$$[G(2,6)]=[ASSOC]+ [\widetilde {ASSOC}]$$ follows from
$\int_{G(2, 6)} \ p_1^2(E)=2$ and $\int_{G(2, 6)} \ p_1(E)e(F)=0$.

{\bf  Theorem 7.5} \ (1) $\frac 12 (p_1(E) + e(F)), \ \frac 12 (p_1(E) - e(F))$ are two generators of
$ H^{4}(G(3,7),{{\bf Z}})$. Their Poincar\' {e} duals are
$[ASSOC]$ and $[\widetilde {ASSOC}]$ respectively;

(2) $\frac 12 (p_1(E)e(F) + e^2(F)), \ \frac 12 (p_1(E)e(F) - e^2(F))\in  H^{8}(G(3,7),{{\bf Z}})$
are  generators and
their Poincar\' {e} duals are
 $[{\bf C}P^2], \ [\overline {{\bf C}P}^2]$ respectively;

(3) $\frac 12 (p_1(E)e(F) + e^2(F)), \ \frac 12 (p_1(E)e(F) - e^2(F))$ and
 $[ASSOC], \ [\widetilde {ASSOC}]$ are dual basis with respect to the universal coefficients Theorem.

The proof is similar to  Theorem 5.5.

\vskip 1cm \noindent{\bf \S 8. \  The case of $G(3,8)$ } \vskip 0.3cm

The Poincar\' {e}  polynomial of Grassmann manifold $G(3,8)$ is
$$p_t(G(3,8)) = (1+t^4+t^8)(1+t^7) = 1+t^4+t^7 + t^8 +t^{11} +t^{15}.$$

Let $E=E(3,8), \ F=F(3,8)$ be canonical vector bundles on $G(3,8)$.
By $\int_{{\bf C}P^2}  \ p_1(E)=1, \ \int_{ASSOC}  \ p_1^2(E)=1$
we know that $$[{\bf C}P^2]\in H_4(G(3,8), {\bf Z}), \ [ASSOC]\in H_8(G(3,8), {\bf Z}),$$
$$p_1(E)\in H^4(G(3,8), {\bf Z}), \
p_1^2(E) \in H^8(G(3,8), {\bf Z})$$ are all generators. By Theorem 3.1, to understand the
structure of the homology and coholomogy groups of dimension $7,11$,  we
need  compute the Poincar\' {e} duals of $p_1(E),  \
p_1^2(E)$.

It is not difficult to compute
$$a=(p_1(E), p_1(E)) = \frac {15}{2 \pi^4} V(G(3,8))=\frac 13 \pi^4.$$
Then $\frac 1a* p_1(E)\in H^{11}(G(3,8), {\bf Z})$ is a generator,
we look for a submanifold $M$ such that
$\int_M \ \frac 1a* p_1(E)=1$.

In \S 7 we define
fibre bundle $\tau\colon\; G(3,8) \to S^7$, $\tau(\pi)=v$ defined by
$\pi A=vA$.  As $v\perp \pi$, $v\pi\in G(4,8)$ and $v\pi A=-A$. Let
 $CAY =\{ \pi\in G(4,8) \ | \ \pi A=-A\}$, called  Cayley submanifold of $G(4,8)$. Then
we have a fibre bundle $\mu\colon\; G(3,8) \to CAY, \ \pi\mapsto v\pi$,
with fibre $S^3=G(3,4)$.
Let  $ASSOC=\{\bar e_1 e_2e_3e_4 \ | \ e_2e_3e_4A=\bar e_1A\}$, where $\bar e_1=(1,0,\cdots,0)$.
Then $M=\mu^{-1}(ASSOC)$ is a $11$ dimensional submanifold of $G(3,8)$.

Let $\bar e_1, e_2,\cdots, e_8$ be $G_2$ frame fields on ${\bf R}^8$,
 $ASSOC$ be generated by $\bar e_1e_2e_3e_4$. Represent
the elements of $\mu^{-1}(\bar e_1e_2e_3e_4)$ by $\tilde e_2\tilde e_3\tilde e_4, \
\tau (\tilde e_2\tilde e_3\tilde e_4)=\tilde e_1$, then
$$\tilde e_1\tilde e_2\tilde e_3\tilde e_4=\bar e_1e_2e_3e_4.$$
Let
$\mbox{d} (\tilde e_2\tilde e_3\tilde e_4) =\sum\limits_{i=2}^4 \ \tilde \omega_i^1\tilde E_{i1} +
\sum\limits_{i=2}^4\sum\limits_{\alpha=5}^8 \ \tilde \omega_i^\alpha\tilde E_{i\alpha}$, \
$\mbox{d}s^2_M = \sum\limits_{i} \ (\tilde \omega_i^1)^2 +
\sum\limits_{i,\alpha} \ (\tilde \omega_i^\alpha)^2$ be the metric on $M$.

Let $\tilde e_1 = \lambda_1\bar e_1+\sum\limits_{i=2}^4 \ \lambda_i e_i, \
\mbox{d} \tilde e_1 = \sum\limits_{i=2}^4 \ \tilde \omega_1^i\tilde e_i +
\sum\limits_{\alpha=5}^8 \ \tilde\omega_1^\alpha e_\alpha$
and $\mbox{d} (
  e_2 e_3 e_4)=\sum\limits_{i,\alpha} \ \omega_i^\alpha E_{i\alpha}$.
It is easy to see    $\tilde\omega_1^\alpha  = \sum \
 \lambda_i\omega_i^\alpha $.
 $\sum\limits_{i=2}^4 \ (\tilde \omega_1^i)^2$ is the metric on the fibres
of $\mu\colon\; M \to ASSOC$ and $\langle \mbox{d} ( e_2 e_3 e_4),\mbox{d} (
  e_2e_3 e_4)\rangle =\sum \ (\omega_i^\alpha)^2$ is metric on $ASSOC$.
From $\langle \bar e_1E_{i\alpha},e_\beta\tilde  e_2\tilde e_3\tilde e_4\rangle \\ =
\lambda_i\delta_{\alpha\beta}$
we have
$$\langle \bar e_1\mbox{d} (e_2e_3e_4),(\mbox{d} \tilde e_1)\tilde
  e_2\tilde e_3\tilde e_4\rangle=
  \langle \sum \ \omega_i^\alpha\bar e_1E_{i\alpha},\sum \ \lambda_j\omega_j^\beta e_\beta
  \tilde  e_2\tilde e_3\tilde e_4\rangle =
\sum\limits_{\alpha} \ (\sum_i \ \lambda_i\omega_i^\alpha)^2 .$$
Hence
\begin{eqnarray*}
\sum\limits \ (\tilde \omega_i^\alpha)^2  &=&
\langle \tilde e_1\mbox{d} (\tilde
  e_2\tilde e_3\tilde e_4),\tilde e_1\mbox{d} (\tilde
  e_2\tilde e_3\tilde e_4)\rangle \\
&=& \langle \bar e_1\mbox{d} (e_2e_3e_4) - (\mbox{d} \tilde e_1)\tilde e_2\tilde e_3\tilde e_4,
\bar e_1\mbox{d} (e_2e_3e_4) - (\mbox{d} \tilde e_1)\tilde e_2\tilde e_3\tilde e_4\rangle\\
&=& \sum_{i,\alpha} \ (\omega_i^\alpha)^2 +
\sum \ (\tilde \omega_1^\alpha)^2 -2 \sum\limits_{\alpha} \ (\sum_i \ \lambda_i\omega_i^\alpha)^2 \\
&=& \sum_{i,\alpha} \ (\omega_i^\alpha)^2 - \sum\limits_{\alpha} \ (\sum_i \ \lambda_i\omega_i^\alpha)^2.\end{eqnarray*}

Then the metric on $M$ can be represented by
$$\mbox{d}s^2_M = \sum_{i,\alpha} \ (\omega_i^\alpha)^2 -
\sum\limits_{\alpha} \ (\sum_i \ \lambda_i\omega_i^\alpha)^2+ \sum\limits_{i=2}^4 \ (\tilde \omega_1^i)^2.$$

For fixed  $\tilde  e_2\tilde e_3\tilde e_4$ of $M$, we can choice $G_2$ frame fields
$\bar e_1, e_2,\cdots, e_8$ such that $\tilde e_1 = \lambda_1\bar e_1+\lambda_2 e_2$. By
$$ \omega_4^5=-\omega_2^7-\omega_3^6, \ \omega_4^6 =\omega_2^8+\omega_3^5, \
 \omega_4^7= \omega_2^5-\omega_3^8,  \ \omega_4^8=-\omega_2^6+\omega_3^7,$$
we have
$$\mbox {d}V_M= \frac 19(1+2\lambda_1^2)^2\tilde \omega_1^2\tilde \omega_1^3\tilde \omega_1^4
\mu^*\mbox {d}V_{ASSOC}.$$
We can show that $\int_{S^3} \ (1+2\lambda_1^2)^2\mbox{d}V_{S^3} = 5\pi^2$, then $V(M)= \frac 23 \pi^6.$

{\bf  Lemma 8.1} \ $\int_M \ \frac 3{\pi^4} *p_1(E)=1$. Then
$[M]$ and $\frac 3{\pi^4} *p_1(E)$ are dual generators of $H_8(G(3,8), {\bf Z}), \ H^8(G(3,8), {\bf Z})$
respectively.

{\bf  Proof } \  By Theorem 3.1, the integration $\int_M \ \frac 1a *p_1(E)$ is an integer. On the other hand
$2\pi^2 p_1(E)$ is a calibration on $G(3,8)$ with comass $\frac 43$, see [7], then we have
 $$\left |\int_M \ \frac 1a *p_1(E)\right |\leq \frac 2{3a\pi^2} V(M)=\frac 43.$$
 We need to show that $\int_M \ \frac 1a *p_1(E)\not =0$.

Let $i\colon\; CAY \to G(4,8)$ be the inclusion, $\mu\colon\; G(3,8)\to CAY$ be sphere bundle associated to
the induced vector bundle $i^*E(4,8) \to CAY$,
$e(i^*E(4,8))=i^*e(E(4,8))\in
H^{4}(CAY, {\bf Z})$ be the Euler class. The induced bundle
$(i\circ\tau)^*E(4,8)\to G(3,8)$ has a nonzero section, $(i\circ\tau)^*e(E(4,8))=0$.
By Gysin sequence for the sphere bundle
$G(3,8)\to CAY$, we can show
$\tau_*(\frac 1a *p_1(E)), (i^*e(E(4,8)))^2$ are two generators of
 $H^8(CAY)$, where $\tau_*\colon\; H^{11}(G(3,8))
 \to H^8(CAY)$ is the integration along the fibre. It is easy to see that
$\int_{ASSOC} \ (i^*e(E(4,8)))^2=0$. If we also have $\int_{ASSOC}
\ \tau_*(\frac 1a *p_1(E))=0$, then $[ASSOC]=0$ in $H_{8}(CAY)$.
This contradict to the fact $i_*[ASSOC]\not =0$.
Then with suitable choice orientation on $M$, we have
$$\int_M \ \frac 1a *p_1(E)= \int_{ASSOC}
\ \tau_*(\frac 1a *p_1(E))=1. \ \ \ \ \Box$$

 Finally we study the $H_7(G(3,8))$ and $H^7(G(3,8))$. Let $I,J,K$ be the quaternion structures on
 ${\bf R}^8={\bf H}^2$ and $Sp(2)$ the symplectic group. As we know  $Sp(2)$ can be viewed as
 a subgroup of $SU(4)\subset SO(8)$, $Sp(2)$ is a subgroup of $Spin_7$. Let $f\colon\; S^7\to G(3,8), \
 f(v)=Iv  Jv Kv$. By $I\bar e_1
 J\bar e_1 K\bar e_1 A =\bar e_2
 \bar e_3 \bar e_4 A =\bar e_1 A$ and $Sp(2)$ acts on $S^7$ transitively,
 we have $ Iv  Jv KvA=vA$ for any $v\in S^7$, \ $\tau(f(v))=v$.

{\bf  Lemma 8.2} \ $b =(p_1^2(E), p_1^2(E))=\frac 52$.

 {\bf  Proof} \ By computation, we have
 \begin{eqnarray*} p_1^2(E) & = & \frac 1{2\pi^4}\{\sum\limits_{i<j}
 \sum\limits_{\alpha <\beta<\gamma<\tau} \ 3\omega_i^{\alpha}
 \omega_j^{\alpha}\omega_i^{\beta}\omega_j^{\beta}\omega_i^{\gamma}
 \omega_j^{\gamma}\omega_i^{\tau}\omega_j^{\tau} \\
&& \qquad + \sum\limits_{j\not = k}
 \sum\limits_{\alpha <\beta,\gamma<\tau} \ \omega_i^{\alpha}
 \omega_j^{\alpha}\omega_i^{\beta}\omega_j^{\beta}\omega_i^{\gamma}
 \omega_k^{\gamma}\omega_i^{\tau}\omega_k^{\tau}\},
 \end{eqnarray*}
 $$i,j,k=2,3,4, \ \alpha,\beta,\gamma,\tau=1,5,6,7,8.$$
 Then
 $$b=(p_1^2(E), p_1^2(E))=\frac 1{4\pi^8}(9\cdot 3\cdot C_5^4 + 3\cdot C_5^2\cdot C_3^2)V(G(3,8))=
 \frac 52. \ \ \ \ \Box$$

 $p_1^2(E)$ is  a harmonic form,
by Theorem 3.1, $\frac 1b *p_1^2(E)\in H^7(G(3,8), {\bf Z})$ is a generator.

{\bf  Lemma 8.3} \ $\frac 25\int_{S^7} \  f^* *p_1^2(E)=1$.  Then
$[f(S^7)]$ and $\frac 25 *p_1^2(E)$ are dual generators of $H_7(G(3,8), {\bf Z}), \ H^7(G(3,8), {\bf Z})$
respectively.

{\bf  Proof} \  Let $e_1, e_2= Ie_1, e_3=Je_1,  e_4= Ke_1, e_5, e_6= Ie_5, e_7=Je_5,  e_8= Ke_5$
be $Sp(2)$ frame fields on ${\bf R}^8$, $f(e_1)=e_2e_3e_4$.
The $1$ forms $\omega_i^\alpha=\langle \mbox{d} e_i, e_\alpha
\rangle, \ i=1,2,3,4, \ \alpha=5,6,7,8$, satisfy
$$\omega_1^{5}=\omega_2^{6}=\omega_3^{7}=\omega_4^{8}, \  \
\omega_1^{6}=-\omega_2^{5}=-\omega_3^{8}=\omega_4^{7},$$
$$ \omega_1^{7}=\omega_2^{8}=-\omega_3^{5}=-\omega_4^{6} , \ \
\omega_1^{8}=-\omega_2^{7}=\omega_3^{6}=-\omega_4^{5}.$$
We have
$$*\omega_2^{5}
 \omega_3^{5}\omega_2^{6}\omega_3^{6}\omega_2^{7}
 \omega_3^{7}\omega_2^{8}\omega_3^{8}|_{f(S^7)}=\omega_2^{1}
 \omega_3^{1}\omega_4^{1}\omega_4^{5}
 \omega_4^{6}\omega_4^{7}\omega_4^{8}|_{f(S^7)}=
\omega^1_{2}
 \omega^1_{3}\omega^1_{4}\omega^1_{5}\omega^1_{6}
 \omega^1_{7}\omega^1_{8}=\mbox{d} V_{S^7}, $$
 $$*\omega_2^{1}
 \omega_3^{1}\omega_2^{6}\omega_3^{6}\omega_2^{7}
 \omega_3^{7}\omega_2^{8}\omega_3^{8}|_{f(S^7)}=-\omega_2^{5}
 \omega_3^{5}\omega_4^{5} \omega_4^{1}
 \omega_4^{6}\omega_4^{7}\omega_4^{8}|_{f(S^7)}=0, $$
$$*\omega_2^{5}
 \omega_3^{5}\omega_2^{6}\omega_3^{6}\omega_2^{7}
 \omega_4^{7}\omega_2^{8}\omega_4^{8}|_{f(S^7)}=\omega_2^{1}
 \omega_3^{1}\omega_4^{1}\omega_4^{5}
 \omega_4^{6}\omega_3^{7}\omega_3^{8}|_{f(S^7)}=\mbox{d} V_{S^7}, $$
 $$*\omega_2^{1}
 \omega_3^{1}\omega_2^{6}\omega_3^{6}\omega_2^{7}
 \omega_4^{7}\omega_2^{8}\omega_4^{8}|_{f(S^7)}=-\omega_2^{5}
 \omega_3^{5}\omega_4^{5} \omega_4^{1}
 \omega_4^{6}\omega_3^{7}\omega_3^{8}|_{f(S^7)}=0, \ \cdots.$$
Then
$$\frac 1b f^**p_1^2(E) = \frac 25\cdot \frac 1{2\pi^4} (3\cdot 3 +3\cdot 2)\mbox{d} V_{S^7}
=\frac 3{\pi^4}\mbox{d} V_{S^7},$$
$$\frac 25\int_{S^7} \  f^* *p_1^2(E)=1.\ \ \ \ \Box$$

Obviously, we also have $\int_{f(S^7)} \ \tau^* \frac 3{\pi^4}\mbox{d} V_{S^7}=1$, then
$$[\tau^* \frac 3{\pi^4}\mbox{d} V_{S^7}]=\frac 25 *p_1^2(E)\in H^7(G(3,8)). $$

{\bf Theorem 8.4} \ \ (1) The Poincar\' {e} dual of $p_1(E)$  is  $[M]$ and
the Poincar\' {e} dual of $\frac 3{\pi^4} *p_1(E)$
  is  $[{\bf C}P^2]$;

(2) The Poincar\' {e} dual of $p_2(F)=p_1^2(E)$  is  $[f(S^7)]$ and
the Poincar\' {e} dual of $\frac 25 *p_1^2(E)$
  is  $[ASSOC]$.

\vskip 1cm \noindent{\bf\S 9. \  The case of $G(4,8)$ } \vskip 0.3cm

Let $E=E(4,8), \ F=F(4,8)$ be canonical vector bundles on Grassmann manifold $G(4,8)$.
We have
$$p_t(G(4,8)) =1+3t^4 +4t^8 + 3t^{12} + t^{16},$$
$$e(E)e(F)=0, \ p_1(E)=-p_1(F), \ p_2(E)= e^2(E), \ p_2(F)= e^2(F),$$
$$p_1^2(E)= p_2(E) + p_2(F), \
p_1(E)p_2(E) = p_1(E)p_2(F)=\frac 12 p_1^3(E),$$ $$ p^2_1(E)e(E) = e^3(E), \ p^2_1(F)e(F) = e^3(F).$$
 By the method used in \S 4,
we can show
$e(E\otimes F)= 6e^4(E)= 6 e^4(F)$. Then by $\int_{G(4,8)} \ e(E\otimes F)=\chi (G(4,8))=12$, we have
$$\int_{G(4,8)} \ e^4(E)= \int_{G(4,8)} \ e^4(F) =2.$$

We first study $4$ and $12$ dimensional cases.
As we know ${\bf C}P^2, \overline {{\bf C} CP}^2, G(2,4)$ can be imbedded in $G(4,8)$ as submanifolds.
Under the map $*\colon\; G(4,8)\to G(4,8)$,  $*{\bf C}P^2$ is  a submanifold of $G(4,8)$.
The following table computes the integration of the characteristic classes on the submanifolds of $G(4,8)$.
\vskip 0.1cm $$\begin{tabular}{|c|c|c|c|}\hline
            &\ \ \ \ ${\bf C}P^2$\ \ \ \ &\ \ \ \ $*{\bf C}P^2$ \ \ &  \ \ $G(2,4)$ \ \ \\ \hline
\ \ \ \ $e(E)$\ \ \ \ &    0       & 1    \  & 0   \\   \hline
        $ e(F)$ &     1        &       0    \ & 0 \\  \hline
       $ p_1(E)$ &     1       &       $-1$    \ & 2 \\  \hline
\end{tabular}$$
 \vskip 0.2cm

Note that $\det \left( \begin{array} {ccccc}  0 & 1 & 0 \\
1 & 0 & 0 \\ 1 & -1 & 2
\end{array} \right) =-2$, as  proof of Theorem 5.5, we can show
$e(E), e(F), p_1(E)\in H^4(G(4,8),{\bf Z})$ or
${\bf C}P^2, *{\bf C}P^2,G(2,4)\in H_4(G(4,8),{\bf Z})$  are  generators.

By Proposition 2.2, $V(G(4,8))= \frac 8{135}\pi^8$ and we can compute
$$a=(e(E),e(E))=(e(F),e(F))=\frac 12(p_1(E),p_1(E))=\frac 4{15}\pi^4,$$
$$(e(E),e(F))=(e(E),p_1(E))=(e(F),p_1(E))=0.$$
 In last section we have Cayley submanifold $CAY =\{ \pi\in G(4,8) \ | \ \pi A=-A\}$
 of $G(4,8)$. The Lie group $Spin_7$ acts on $CAY$ transitively.
Let $ e_1, e_2,\cdots, e_8$ be $Spin_7$ frame fields on ${\bf R}^8$,
then $CAY$ be generated by $e_1e_2e_3e_4$.  By the equations listed in the proof of Lemma 7.3, we have
\begin{eqnarray*}  && \mbox {d}(e_1e_2e_3e_4) = \sum\omega_i^\alpha E_{i\alpha} \\
 & & =\omega^5_1(E_{15}+E_{48})+\omega^6_1(E_{16}+E_{47})+\omega^7_1(E_{17}-E_{46})
   +\omega^8_1(E_{18}-E_{45}) \\ && \quad
 +\omega^5_2(E_{25}+E_{47})+\omega^6_2(E_{26}-E_{48})
 +\omega^7_2(E_{27}-E_{45})+\omega^8_2(E_{28}+E_{46}) \\
&& \quad +\omega^5_3(E_{35}+E_{46})+\omega^6_3(E_{36}-E_{45})+
\omega^7_3(E_{37}+E_{48})+\omega^8_3(E_{38}-E_{47}).
 \end{eqnarray*}
Then the induced metric is
\begin{eqnarray*}\mbox {d}s^2_{CAY} & = &
(\omega^5_1-\omega^6_2)^2+(\omega^5_1+\omega^7_3)^2+(\omega^6_2-\omega^7_3)^2\\
&& +(\omega^6_1+\omega^5_2)^2+(\omega^6_1-\omega^8_3)^2
 +(\omega^5_2-\omega^8_3)^2 \\
&& +(\omega^7_1-\omega^5_3)^2+(\omega^7_1-\omega^8_2)^2 +(\omega^8_2+\omega^5_3)^2 \\
&& +
(\omega^8_1+\omega^7_2)^2+(\omega^8_1+\omega^6_3)^2
 +(\omega^7_2+\omega^6_3)^2. \end{eqnarray*}
By $(\omega^5_1-\omega^6_2)(\omega^5_1+\omega^7_3)(\omega^6_2-\omega^7_3)=
2\omega^5_1\omega^7_3\omega^6_2, \ \cdots$, we get
$$\mbox {d}V_{CAY}
 =16\omega^5_1\omega^5_2\omega^5_3\ \cdots\ \omega^8_1\omega^8_2\omega^8_3.$$

As shown in [7], $2\pi^2 p_1(E)$ is a calibration on $G(4,8)$ with comass $\frac 32$ and
$CAY$ is a calibrated submanifold of $2\pi^2 * p_1(E)$, then
$$2\pi^2 * p_1(E)|_{CAY} = \frac 32 \mbox {d}V_{CAY}.$$
By triality transformation, we can show that $CAY$ is isometric to $G(3,7)$, then
$V(CAY)=V(G(3,7))=\frac {16\pi^6}{45}$.
These shows
$$\int_{CAY} \ \frac 1{2a} * p_1(E) =  \int_{CAY} \ \frac {45}{32\pi^6} \mbox {d}V_{CAY}=\frac 12.$$
By Theorem 3.1,
$e(E), e(F), p_1(E)$ can not be generators of $H^4(G(4,8),{\bf Z})$.
We have proved

{\bf Lemma 9.1} \  ${\bf C}P^2, *{\bf C}P^2, G(2,4)$ are the generators of  $H_4(G(4,8),{\bf Z})$,
the dual generators are the first column in the following table.
\vskip 0.1cm $$\begin{tabular}{|c|c|c|c|}\hline
            &\ \ \ \ ${\bf C}P^2$\ \ \ \ &\ \ \ \ $*{\bf C}P^2$ \ \ &  \ \ $G(2,4)$ \ \ \\ \hline
\ \ \ \ $e(E)$\ \ \ \ &    0       & 1    \  & 0   \\   \hline
        $ e(F)$ &     1        &       0    \ & 0 \\  \hline
       $\frac 12(p_1(E)+e(E)-e(F))$ &     0       &       $0$    \ & 1 \\  \hline
\end{tabular}$$
 \vskip 0.2cm

The inner product of $e(E),  e(F),  \frac 12(p_1(E) +e(E) -e(F))$ form a matrix
$$A=\frac  {4\pi^4}{15} \left( \begin{array} {ccc}  1 & 0 & \frac 12 \\
0 & 1 & -\frac 12 \\ \frac 12 & -\frac 12 & 1
\end{array} \right), \  \ A^{-1} = \frac  {15}{4\pi^4} \left( \begin{array} {ccc}
 \frac 32 & -\frac 12 & -1 \\
-\frac 12 & \frac 32 & 1 \\ -1 & 1 & 2
\end{array} \right).  $$

{\bf Lemma 9.2} \ (1) $*e(E) =\frac  {2\pi^4}{15}e^3(E), \ *e(F) =\frac  {2\pi^4}{15}e^3(F), \
*p_1(E) =\frac  {2\pi^4}{15}p_1^3(E);$

(2) $\frac 12( e^3(E) -\frac 12 p_1^3(E)),
\frac 12( e^3(F) +\frac 12 p_1^3(E)),
\frac 12 p_1^3(E)\in H^{12}
(G(4,8), {\bf Z})$ are the generators.

{\bf Proof} \  $*e(E), *e(F), *p_1(E)$ and $e^3(E),e^3(F), p_1^3(E)$ are two
generators of the coholomogy group $H^{12}(G(4,8))$, they are all the harmonic forms on $G(4,8)$.
Then  $e^3(E),e^3(F), p_1^3(E)$ can be represented by  $*e(E), *e(F), *p_1(E)$. Assuming $e^3(E) =
\lambda *e(E) +\mu *e(F) +\nu *p_1(E)$, by $e(F)\wedge *e(E)=0, \ e(F)\wedge *p_1(E)=0$,
we have $\mu=\nu=0$,
$$2=\int_{G(4,8)} \ e(E)\wedge e^3(E) =\lambda \int_{G(4,8)} \  e(E)\wedge *e(F) = \lambda\frac  {4\pi^4}{15}.$$
Then $*e(E) =\frac  {2\pi^4}{15}e^3(E)$. The other two equalities can be proved similarly.

Then
\begin{eqnarray*} && (*e(E), *e(F),*\frac 12(p_1(E) +e(E) -e(F)))
A^{-1} \\
&& = (\frac 12( e^3(E) -\frac 12 p_1^3(E)), \frac 12( e^3(F) +\frac 12 p_1^3(E)),
\frac 12 p_1^3(E)). \ \ \ \ \Box\end{eqnarray*}

{\bf Lemma 9.3} \ The following table compute the integration of the characteristic classes on the submanifolds of $G(4,8)$.
\vskip 0.1cm $$\begin{tabular}{|c|c|c|c|}\hline
            &\ \ \ \ $G(4,7)$\ \ \ \ &\ \ \ \ $G(3,7)$ \ \ &  \ \ $CAY$ \ \ \\ \hline
\ \ \ \ $e^3(E)$\ \ \ \ &    2      & 0    \  & $-1$   \\   \hline
        $ e^3(F)$ &     0        &       2    \ & 1 \\  \hline
       $ p_1^3(E)$ &     0       &       $0$    \ & 2 \\  \hline
\end{tabular}$$
 \vskip 0.2cm

{\bf Proof} \ The first column follows from
$\int_{G(4,7)} \ e^{3}(E(4,8))=\int_{G(4,7)} \ e^{3}(E(4,7))=2$
and $ e(F(4,8))|_{G(4,7)}=0, \  p_1^3(E)= 2p_1(E)p_2(F)=2p_1(E)e^2(F)$.  The second column can be
proved similarly. For third column, we have proved
$\int_{CAY} \ \frac 1{2a} * p_1(E) =  \frac 12$,
then $\int_{CAY} \  p^3_1(E) =  2$.
By compute $*e(E)|_{CAY}, \ *e(F)|_{CAY}$, we can show
 $\int_{CAY} \  e^3(E) =  -1$
and $\int_{CAY} \  e^3(F) = 1$. \ \ \ \ $\Box$

{\bf Theorem 9.4} \ (1) $e(E),  e(F),  \frac 12(p_1(E) +e(E) -e(F))\in H^4(G(4,8), {\bf Z})$
are the generators, their dual generators  are
$[{\bf C}P^2], [*{\bf C}P^2], [G(2,4)]\in H_4(G(4,8), {\bf Z})$;

(2) $\frac 12e^3(E), \frac 12 e^3(F),
\frac 12 p_1^3(E)$ and $[G(4,7)],  [G(3,7)], [CAY]$ are the generators of  $H^{12}
(G(4,8), {\bf Z})$ and $H_{12}
(G(4,8), {\bf Z})$  respectively;

(3) The Poincar\' {e} duals of $e(E),  e(F),  \frac 12(p_1(E) +e(E) -e(F))$ are
$$[G(4,7)], \ [G(3,7)], \ [CAY]+[G(4,7)]-[G(3,7)]$$ respectively.

{\bf Proof} \  By Lemma 9.2,  $\frac 12( e^3(E) -\frac 12 p_1^3(E)),
\frac 12( e^3(F) +\frac 12 p_1^3(E)),
\frac 12 p_1^3(E)$ are the generators of $H^{12}
(G(4,8), {\bf Z})$. Then $\frac 12e^3(E),
\frac 12 e^3(F),
\frac 12 p_1^3(E)$ are also the generators of $H^{12}
(G(4,8), {\bf Z})$. \ \ \ \ $\Box$

By Theorem 3.1, we can compute the  Poincar\' {e} duals of
$$\frac 12( e^3(E) -\frac 12 p_1^3(E)), \frac 12( e^3(F) +\frac 12 p_1^3(E)),
\frac 12 p_1^3(E).$$

By Theorem 9.4, $\frac 12(p_1(E) +e(E) -e(F))e(E)=
\frac 12(p_1(E)e(E) + e^2(E))$ and $\frac 12(p_1(E)e(F) - e^2(F)), \
\frac 12(p_1(E)e(E) - e^2(E)), \ \frac 12(p_1(E)e(F) + e^2(F))$
are integral cocycles.
The submanifolds $ASSOC, \ \widetilde {ASSOC}$  defined in \S 7
are also the submanifolds of $G(4,8)$, then $* ASSOC, \ * \widetilde {ASSOC}$
are  submanifolds of $G(4,8)$. The following table can be proved by  Lemma 7.4.

 \vskip 0.1cm $$\begin{tabular}{|c|c|c|c|c|}\hline
 & \  $ASSOC$ \  & \   $\widetilde {ASSOC}$ \ \ &   \ $*ASSOC$ \
 & \ $* {\widetilde {ASSOC}}$ \ \\ \hline
\  $\frac 12(e^2(F)+ p_1(E)e(F) )$  \ &    1      & 0    \  & 0 & 0   \\   \hline
 $\frac 12(e^2(F) - p_1(E)e(F))$ &     0        &      $1$    \ & 0 &  0  \\  \hline
$\frac 12(e^2(E) + p_1(E)e(E))$ &     0       &       $0$    \ & 1 & 0 \  \\  \hline
$\frac 12(e^2(E) - p_1(E)e(E))$ &     0       &       $0$    \ & $0$ & $1$ \\  \hline
\end{tabular}$$
 \vskip 0.2cm

{\bf Theorem 9.5} \ The characteristic classes
$$\frac 12(e^2(F)+  p_1(E)e(F)),\ \frac 12(e^2(F)-  p_1(E)e(F)),$$
$$ \frac 12(e^2(E) + p_1(E)e(E)), \  \frac 12(e^2(E) - p_1(E)e(E))$$ are the
 generators of $H^{8}
(G(4,8), {\bf Z})$. Their Poincar\' {e} duals are
 $$[ASSOC], \ [\widetilde {ASSOC}], \ [*ASSOC], \ [* {\widetilde {ASSOC}}]$$ respectively.

{\bf Proof} \  To see the Poincar\' {e} dual of
$\xi = \frac 12(e^2(F)+ p_1(E)e(F))$ is $ASSOC$, we want to show that for any $\eta\in
H^{8}(G(4,8))$ we have $\int_{G(4,8)} \ \xi\wedge \eta =\int_{ASSOC} \ \eta$.
We can take $\eta =\frac 12(e^2(F)\pm p_1(E)e(F)), \ \frac 12(e^2(E) \pm p_1(E)e(E))$
 to verify this equation. \ \ \ \ $\Box$

By ${\bf R}^8 = {\bf R}^3\oplus{\bf R}^5$, we see the product Grassmann $G(2,3)\times
G(2,5), G(1,3)\times G(3,5)$ can imbedded in $G(4,8)$ and we have
\vskip 0.1cm $$\begin{tabular}{|c|c|c|c|c|}\hline
 & \  $G(4,6)$ \  & \   $G(2,6)$ \ \ &   \ $G(2,3)\times
G(2,5)$ \
 & \ $G(1,3)\times G(3,5)$ \ \\ \hline
\  $e^2(E)$  \ &    $2$      & $0$      & 0 & 0   \\   \hline
 $e^2(F)$ &     $0$        &      $2$     & 0 &  0  \\  \hline
 \ $ p_1(E)e(E)$ \ &     0       &       $0$     & $4$ & 0   \\  \hline
$ p_1(E)e(F)$ &     0       &       $0$     & $0$ & $4$ \\  \hline
\end{tabular}$$
 \vskip 0.2cm

Then  $$G(4,6), \ G(2,6), \ G(2,3)\times
G(2,5), \ G(1,3)\times G(3,5) \in H_8(G(4,8),{\bf R})$$
 and $$e^2(E), \ e^2(F), \ p_1(E)e(E), \ p_1(E)e(F)
 \in H^8(G(4,8),{\bf R})$$ are
also the generators.

As application, we consider the immersion $f\colon\; M \to {{\bf R}}^8$ of a compact oriented $4$
dimensional manifold, $g\colon\; M \to G(4,8)$ be its Gauss map.
From following table
\vskip 0.1cm $$\begin{tabular}{|c|c|c|c|}\hline
            &\ \ \ \ $G(4,5)$\ \ \ \ &\ \ \ \ $G(1,5)$ \ \ &  \ \ $G(2,4)$ \ \ \\ \hline
\ \ \ \ $e(E)$\ \ \ \ &    2       & 0    \  & 0   \\   \hline
        $ e(F)$ &     0        &      2    \ & 0 \\  \hline
       $ p_1(E)$ &     0       &       0    \ & 2 \\  \hline
\end{tabular}$$
 \vskip 0.2cm
\noindent we have
$$g_*[M] =\frac 12\chi(M)[G(4,5)] + \lambda [G(1,5)] + \frac 32 Sign(M)[G(2,4)],$$
where $\lambda = \frac 12\int_{M} \ e(F(4,8))$ and $Sign(M)$ be the Signature of $M$. $\lambda=0$
if $f$ is an imbedding.

If $g$ is the Gauss map of immersion of in $ {{\bf R}}^7$ or $ {{\bf R}}^6$,
we have
$$g_*[M] =\frac 12\chi(M)[G(4,5)]  + \frac 32 Sign(M)[G(2,4)].$$

\vskip 1cm \noindent{\bf \S 10. \  The cohomology groups on $ASSOC$ } \vskip 0.3cm

The submanifold $ASSOC\approx G_2/SO(4)$ of Grassmann manifold $G(3,7)$ is
important in theory of calibration, see  [7], [9].
In [4]  Borel and Hirzebruch studied the characteristic classes on homogenous spaces,  they
computed the cohomology of $ASSOC$.
In the following we use Gysin sequence to study the cohomology of $ASSOC$.

 As \S 7, let $G(2,7)$ and  $G(3,7)$ be Grassmann manifolds
on ${\bf R}^7\subset {\bf R}^8$ generated by
$\bar e_2, \cdots,   \bar e_8$, and $S^6\subset {\bf R}^7$ the unit sphere.
 There is a fibre bundle
$\tau_1\colon\; G(2,7) \to S^6$ defined by $\pi A=\bar e_1vA, \ \tau_1(\pi)=v$, where $A\in C\ell_8$
is defined in \S 7.
 For any $G\in G_2$,
we have the following commutative diagram
$$\begin{array} {cccc}   G(2, 7) &  \stackrel{G}\longrightarrow & G(2,   7)
\\ \tau_1 \downarrow
& & \downarrow \tau_1 \\ S^6 &  \stackrel{G} \longrightarrow &
S^6. \end{array}$$
The fibre $\tau_1^{-1}(\bar e_2) =\{ v\wedge Jv \ | \ v\in S^6, \ v\perp \bar e_2\}\approx
{\bf C}P^2$, see [16].

 Then for any $\pi\in G(2,7), \ v=\tau_1(\pi),
\  v\wedge \pi\in ASSOC$. This defines map
$$\tau_2\colon\; G(2,7) \to ASSOC, \ \ \pi\mapsto v\wedge \pi. $$
For any $e_2e_3e_4\in ASSOC, \ e_2e_3e_4 A=\bar e_1 A$,
then $\tau_2(e_3e_4)=e_2e_3e_4$. This shows

{\bf Lemma 10.1} \ $\tau_2\colon\; G(2,7) \to ASSOC$ is a fibre bundle with fibre $G(2,3)=S^2$.

Let $i\colon\; ASSOC \to G(3,7)$ be inclusion. It is easy to see
$G(2,7)$  isomorphic to the sphere bundle $S(\tilde E) =\{ v\in \tilde E, \ |v|=1\}$ of
the induced bundle $\tilde E=i^*E(3,7)$.
Let $e(E(3,7))\in H^3(G(3,7), {\bf Z})$ be the Euler class of $E(3,7)$,
$2e(E(3,7))=0$, see [12] p.95-103. Then $e(\tilde E)=i^*e(E(3,7))\in H^3(ASSOC,  {\bf Z})$ is the
 Euler class of the induced bundle $\tilde E$.
There is a Gysin exact sequence for the sphere bundle $G(2, 7) \to ASSOC$,
\begin{eqnarray*} &&  \longrightarrow H^q(ASSOC) \stackrel{\tau_2^* }
\longrightarrow  H^q(G(2, 7)) \stackrel{\tau_{2*}}\longrightarrow H^{q-2}(ASSOC) \\
&& \quad  \stackrel{\wedge e(\tilde E)}  \longrightarrow H^{q+1}(ASSOC)\stackrel{\tau_2^* }
\longrightarrow  H^{q+1}(G(2, 7))\longrightarrow, \end{eqnarray*}
 where $\tau_{2*}$ is the integration along the fibre.
The coefficients of the cohomology groups can be ${\bf R}, {\bf Z}$ or ${\bf Z}_2$.

{\bf Lemma 10.2} \ $e(\tilde E)=i^*e(E(3,7))\not = 0$.

{\bf Proof} \ The map $\tau_{2*}\colon\;  H^q(G(2,7),  {\bf Z}) \to  H^{q-2}(ASSOC,  {\bf Z})$ is the
integration along the fibre. Let $\bar e_1,
e_2, e_3, \cdots, e_8$ be $G_2$ frame fields, $G(2,7)$ is generated by $e_3e_4$ and $\tau_2(e_3e_4)=
e_2e_3e_4$.  Then the Euler class of vector bundle $E(2,7)$ can be represented by
 $$e(E(2,7)) =\frac 1{2\pi} \omega_3^2\wedge \omega_4^2 + \frac 1{2\pi} \sum\limits_{\alpha =5}^8 \
\omega_3^\alpha\wedge\omega_4^\alpha$$
and $\omega_3^2\wedge\omega_4^2$ is the volume element of the fibre at $e_3e_4$. Then
$\tau_{2*}(e(E(2,7))=2$.

By Gysin sequence, the map
$\tau_{2*}\colon\; H^2(G(2, 7)) \to H^{0}(ASSOC)$ is  surjective if $e(\tilde E)= 0$.
These contradicts to the fact  $\tau_{2*}(e(E(2,7))=2$ and $e(E(2,7))\in H^2(G(2, 7), {\bf Z})$
be a generator. \ \ \ \ $\Box$

Then $e(E(3,7))$ is also nonzero and a torsion element of $H^3(G(3,7),  {\bf Z})$.

{\bf Theorem 10.3} \ The cohomology groups of  $ASSOC$ are
$$H^q(ASSOC,  {\bf Z}_2)= \left\{ \begin{array} {ll} {\bf Z}_2, \ \ &  q\not= 1,7, \\
0, \ &  q= 1,7; \end{array} \right.$$
$$H^q(ASSOC,  {\bf Z})= \left\{ \begin{array} {ll} {\bf Z}, \ \ &  q= 0,4,8, \\
{\bf Z}_2, \ &  q= 3,6, \\ 0, \ & q= 1,2,5,7; \end{array} \right.$$
$$H^q(ASSOC,  {\bf R})= \left\{ \begin{array} {ll} {\bf R}, \ \ &  q= 0,4,8, \\
0, \ &  q\not = 0,4,8. \end{array} \right.$$

{\bf Proof} \  $G(2,7)$ is a Kaehler  manifold, the coholomogy of $G(2,7)$ is
generated by Euler class $e(E(2,7))$. We prove the case of  ${\bf Z}_2$ coefficients,
the other cases are left to the reader.

By Gysin sequence, we have
$$  0= H^{-2}(ASSOC)\stackrel{\wedge e(\tilde E)}  \longrightarrow H^{1}(ASSOC)\stackrel{\tau_2^* }
\longrightarrow  H^{1}(G(2, 7))=0,$$
\begin{eqnarray*}  && 0= H^{-1}(ASSOC) \stackrel{\wedge e(\tilde E)}\longrightarrow H^{2}(ASSOC)\stackrel{\tau_2^* }
\longrightarrow  H^{2}(G(2, 7))\stackrel{\tau_{2*}=0 }
\longrightarrow  H^0(ASSOC) \\
&&\quad \stackrel{\wedge e(\tilde E)}  \longrightarrow H^{3}(ASSOC)
\stackrel{\tau_2^* }
\longrightarrow  H^{3}(G(2, 7))=0.\end{eqnarray*}
These shows $H^{1}(ASSOC)=0$ and $H^{2}(ASSOC)\cong H^{2}(G(2, 7)), \
H^0(ASSOC)\cong H^{3}(ASSOC)$.
By
$$ 0= H^{1}(ASSOC) \stackrel{\wedge e(\tilde E)}\longrightarrow H^{4}(ASSOC)\stackrel{\tau_2^* }
\longrightarrow  H^{4}(G(2, 7))\stackrel{\tau_{2*}}
\longrightarrow  H^2(ASSOC) $$
and $\tau_{2*}=0\colon\; H^{4}(G(2, 7), {\bf Z}_2)\to H^{4}(ASSOC, {\bf Z}_2)$,
we have $$ H^{4}(ASSOC)\cong H^{4}(G(2, 7)).$$
The cases of $q=5,\cdots,8$ can be proved similarly or by  Poincar\'{e} duality.

\vskip 1cm \centerline {\large References} \vskip 0.3cm

{\small

\noindent [1] \ Chen, W.H.,    \ The differential geometry of
Grassmann manifold as submanifold, Acta Math., Sinica, {\bf 31A}(1988),
46-53.

\noindent [2] \  Chern, S.S.,  Characteristic classes of Hermitian manifolds.
Annals of. Math.  {\bf 47}(1946), 85-121.

\noindent [3] \  Chern, S.S. and Spanier, E.H.,  A theorem on orientable surfaces in
four-dimensional space.
Comm. Math. Helv. {\bf 25}(1951), 205-209.

\noindent [4] \  Borel, A. and Hirzebruch, F., Characteristic classes and homogenous spaces, I. Amer. J. Math.
{\bf 80}(1958), 458-538.

\noindent [5] \  Bott, R. and  Tu, L.W.,  \  Differential forms in algebraic topology,
Springer-Verlag, New York, GTM {\bf 82},  1982.

\noindent [6] \  Bryant, R.L.,  \   Submanifolds and special
structures on the octonians,   Diff Geom.  J.,  {\bf 17}(1982),
185-232.

\noindent [7] \  Gluck, H.,  Mackenzie, D. and  Morgan, F.,
   Volume-minimizing cycles in Grassmann  manifolds,    Duke
Math.  J.,   {\bf 79}(1995),   335-404.

\noindent [8] \  Greub, W.,  Halperin, S.  and Vanstone, R., \
Connections, curvature, and cohomology, Academic Press, New Yark, 1976.

\noindent [9] \  Harvey, F.R.,    and  Lawson,  Jr H.B.,
Calibrated Geometries,   Acta Math. J., {\bf 148}(1982),   47-157.

\noindent [10] \ Lawson, Jr  H.B., Michelsohn,  M.,  Spin
geometry. Princeton University Press. Princeton New Jersey, 1989

\noindent [11]  \ Kabayashi, S. and Nomizu, K., \ Foundations of
differential geometry, vol. 2,  Interscience Publishers, New York,
1969.

\noindent [12]  \
 Milnor, J.W. and Stasheff, J.D., Characteristic classes, Ann. Math. Studies., No. {\bf 76},
Princeton, 1974.

\noindent [13]  \ Wolf, J.A., Spaces of constant curvature, McGraw-Hill, New Yark, 1967.

\noindent [14] \    Zhou, J.W., \ Irreducible Clifford Modules,
Tsukuba J. Math., {\bf 27}(2003), 57-75.

\noindent [15] \   Zhou, J.W., \ Spinors, Calibrations and
Grassmannians, Tsukuba J. Math.,  {\bf 27}(2003), 77-97.

\noindent [16] \ Zhou, J.W. and Huang, H., Geometry on Grassmann Manifolds G(2,8) and G(3,8),
Math. J. of Okayama Univ., {\bf 44}(2002), 171-179.

\noindent [17] \  Zhou, J.W.,  \  Morse functions on Grassmann manifolds, Proc. of the Royal Soc. of Edinburgh,
{\bf 135A}(2005), 1-13.

\noindent [18] \ Zhou, J.W., \  Totally geodesic submanifolds in Lie groups, Acta Math. Hungarica.£¬
{\bf 111}(2006), 29-41.

\noindent [19] \ Zhou, J.W., \  The geometry and topology on Grassmann manifolds, Math. J. Okayama Univ,.
{\bf 48}(2006), 181-195.

\noindent [20] \ Zhou, J.W., \  A note on characteristic classes, Czechoslovak Math. J., 56£¨2006£©£¬721-732.

 }

\end{document}